\documentclass[acmtog,copyrightmode]{acmart}

\usepackage{booktabs} 

\usepackage{subfig}
\usepackage{multirow}
\usepackage{wrapfig}
\usepackage{algorithmic}
\usepackage{amsfonts}
\usepackage{amsmath}

\usepackage[ruled]{algorithm2e} 

\SetAlFnt{\small}
\SetAlCapFnt{\small}
\SetAlCapNameFnt{\small}
\SetAlCapHSkip{0pt}

\acmJournal{TOG} 

\hypersetup{draft}
\begin{document}

\title{Evaluations of The Hierarchical Subspace Iteration Method}

\author{Ahmad Nasikun}
\affiliation{%
  \institution{Delft University of Technology}
  \department{Department of Intelligent Systems}
  \city{Delft}
  \country{The Netherlands}
}
\affiliation{%
  \institution{Universitas Gadjah Mada}
   \department{Department of Electrical and Information Engineering}
  \city{Yogyakarta}
  \country{Indonesia}}
\email{ahmad.nasikun@ugm.ac.id}
\author{Klaus Hildebrandt}
\affiliation{%
  \institution{Delft University of Technology}
  \department{Department of Intelligent Systems}
  \city{Delft}
  \country{The Netherlands}
}
\email{K.A.Hildebrandt@tudelft.nl}
\renewcommand\shortauthors{A. Nasikun \& K. Hildebrandt}

\settopmatter{printacmref=false}

\begin{abstract}
\textit{This document contains additional experiments concerned with the evaluation of the Hierarchical Subspace Iteration Method, which is introduced in~\cite{Nasikun2021}}.  
\end{abstract}

\maketitle

\newcommand{\red}[1] {{\color{red}{{#1}}}}
\newcommand{\eg}{\textit{e.g. }}
\newcommand{\ie}{\textit{i.e. }}
\newcommand{\etal}{\textit{et al. }}
\definecolor{minorRevColor}{rgb}{0.,0.,0.}
\newcommand{\minorrev}[1]{{\color{minorRevColor}{#1}}}

\newcommand{\level}[4]{{#1}^{#2}_{{#3}_{#4}}}		
\newcommand{\refFigErrorTolerance}{1}

\newcommand{\refFigErrorAnalytic}{7}
\newcommand{\refFigErrFlip}{8}

\section{Justification of Design Choices}
In this section, we present experiments that address the justification of our design and evaluation choices for the Hierarchical Subspace Iteration Method (HSIM).
\begin{table}[t]
\begin{tabular}{|l|l|r|r|r|r|r|}
\hline
\multicolumn{1}{|c|}{\multirow{2}{*}{\shortstack[l]{Model\\(\#Verts.)}}} & \multicolumn{1}{c|}{\multirow{2}{*}{Type}} & \multicolumn{1}{c|}{\multirow{2}{*}{Acc.}} & \multicolumn{1}{c|}{\multirow{2}{*}{Hier.}} & \multicolumn{2}{c|}{Solve}                            & \multicolumn{1}{c|}{\multirow{2}{*}{Total}} \\ \cline{5-6}
\multicolumn{1}{|c|}{}                       & \multicolumn{1}{c|}{}                      & \multicolumn{1}{c|}{}                      & \multicolumn{1}{c|}{}                       & \multicolumn{1}{c|}{Iter} & \multicolumn{1}{c|}{Time} & \multicolumn{1}{c|}{}                       \\ \hline \hline
\multirow{6}{*}{\shortstack[l]{Gargoyle\\(85k)}}              & \multirow{2}{*}{Dijkstra}                  & 1e-2                                       & \multirow{2}{*}{2.2}                        & F|1|1                     & 8.4                       & 10.5                                        \\ \cline{3-3} \cline{5-7} 
                                             &                                            & 1e-4                                       &                                             & F|3|3                     & 15.8                      & 18.0                                        \\ \cline{2-7} 
                                             & \multirow{2}{*}{STVD}                      & 1e-2                                       & \multirow{2}{*}{4.0}                        & F|1|1                     & 8.9                       & 13.0                                        \\ \cline{3-3} \cline{5-7} 
                                             &                                            & 1e-4                                       &                                             & F|4|4                     & 20.8                      & 24.8                                        \\ \cline{2-7} 
                                             & \multirow{2}{*}{Heat M.}               & 1e-2                                       & \multirow{2}{*}{203.6}                      & F|1|1                     & 10.4                      & 214.0                                       \\ \cline{3-3} \cline{5-7} 
                                             &                                            & 1e-4                                       &                                             & F|4|3                     & 22.6                      & 226.2                                       \\ \hline
\multirow{6}{*}{\shortstack[l]{Fertility\\(270k)}}            & \multirow{2}{*}{Dijkstra}                  & 1e-2                                       & \multirow{2}{*}{9.6}                                         & F|1|1                     & 26.2                      & 35.8                                        \\ \cline{3-3} \cline{5-7} 
                                             &                                            & 1e-4                                       &                                             & F|3|4                     & 59.3                      & 68.9                                        \\ \cline{2-7} 
                                             & \multirow{2}{*}{STVD}                      & 1e-2                                       & \multirow{2}{*}{30.4}                       & F|1|1                     & 27.3                      & 57.7                                        \\ \cline{3-3} \cline{5-7} 
                                             &                                            & 1e-4                                       &                                             & F|3|2                     & 44.8                      & 75.2                                        \\ \cline{2-7} 
                                             & \multirow{2}{*}{Heat M.}               & 1e-2                                       & \multirow{2}{*}{1218.8}                     & F|1|1                     & 30.8                      & 1249.6                                      \\ \cline{3-3} \cline{5-7} 
                                             &                                            & 1e-4                                       &                                             & F|3|3                     & 53.7                      & 1272.5                                      \\ \hline
\multirow{6}{*}{\shortstack[l]{Oil-pump\\(570k)}}              & \multirow{2}{*}{Dijkstra}                  & 1e-2                                       & \multirow{2}{*}{30.3}                       & F|1|1                     & 63.6                      & 93.9                                       \\ \cline{3-3} \cline{5-7} 
                                             &                                            & 1e-4                                       &                                             & F|4|3                     & 112.6                     & 142.8                                       \\ \cline{2-7} 
                                             & \multirow{2}{*}{STVD}                      & 1e-2                                       & \multirow{2}{*}{100.1}                      & F|1|1                     & 65.6                     & 166.6                                       \\ \cline{3-3} \cline{5-7} 
                                             &                                            & 1e-4                                       &                                             & F|4|3                     & 113.1                     & 214.1                                       \\ \cline{2-7} 
                                             & \multirow{2}{*}{Heat M.}               & 1e-2                                       & \multirow{2}{*}{4339.4}                           & F|2|2                           & 88.9                          & 4628.3                                            \\ \cline{3-3} \cline{5-7} 
                                             &                                            & 1e-4                                       &                                         & F|4|4                            &  136.6                           & 4676.0                                             \\ \hline
\end{tabular}
\caption{The timings and iteration counts for computing 100 eigenpairs on different meshes with three different schemes for approximating the geodesic distances are shown. The timings for the construction of the hierarchy are additionally listed.  
	}
	\label{tab:basisType}
\end{table}

\begin{figure}[]
    \includegraphics[width=\linewidth]{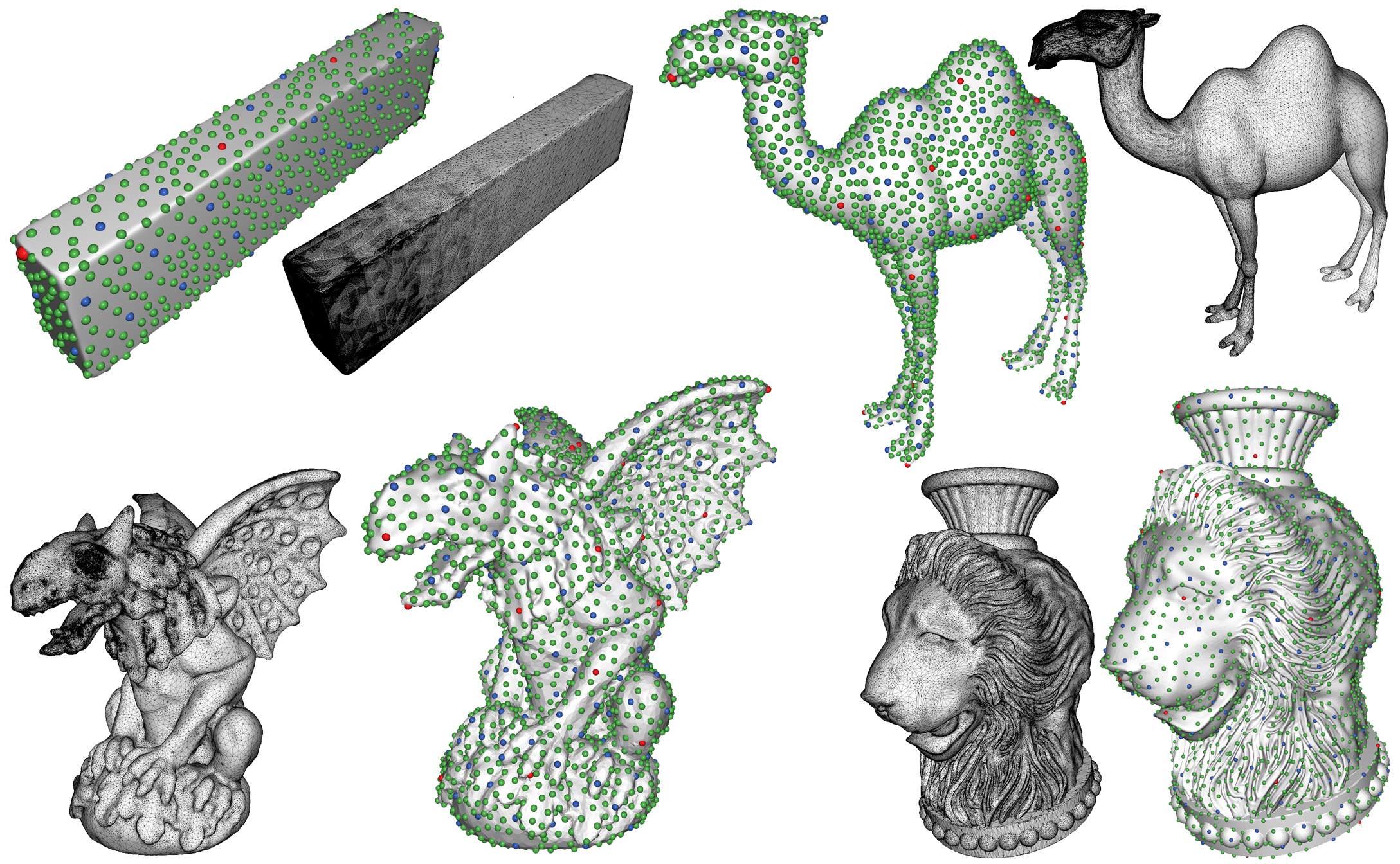}
    \caption{Results of vertex hierarchy construction using farthest point sampling on meshes with spatially varying sampling densities are shown.}
    \label{fig:sampling} 
\end{figure}

\begin{figure}[]
    \includegraphics[width=\linewidth]{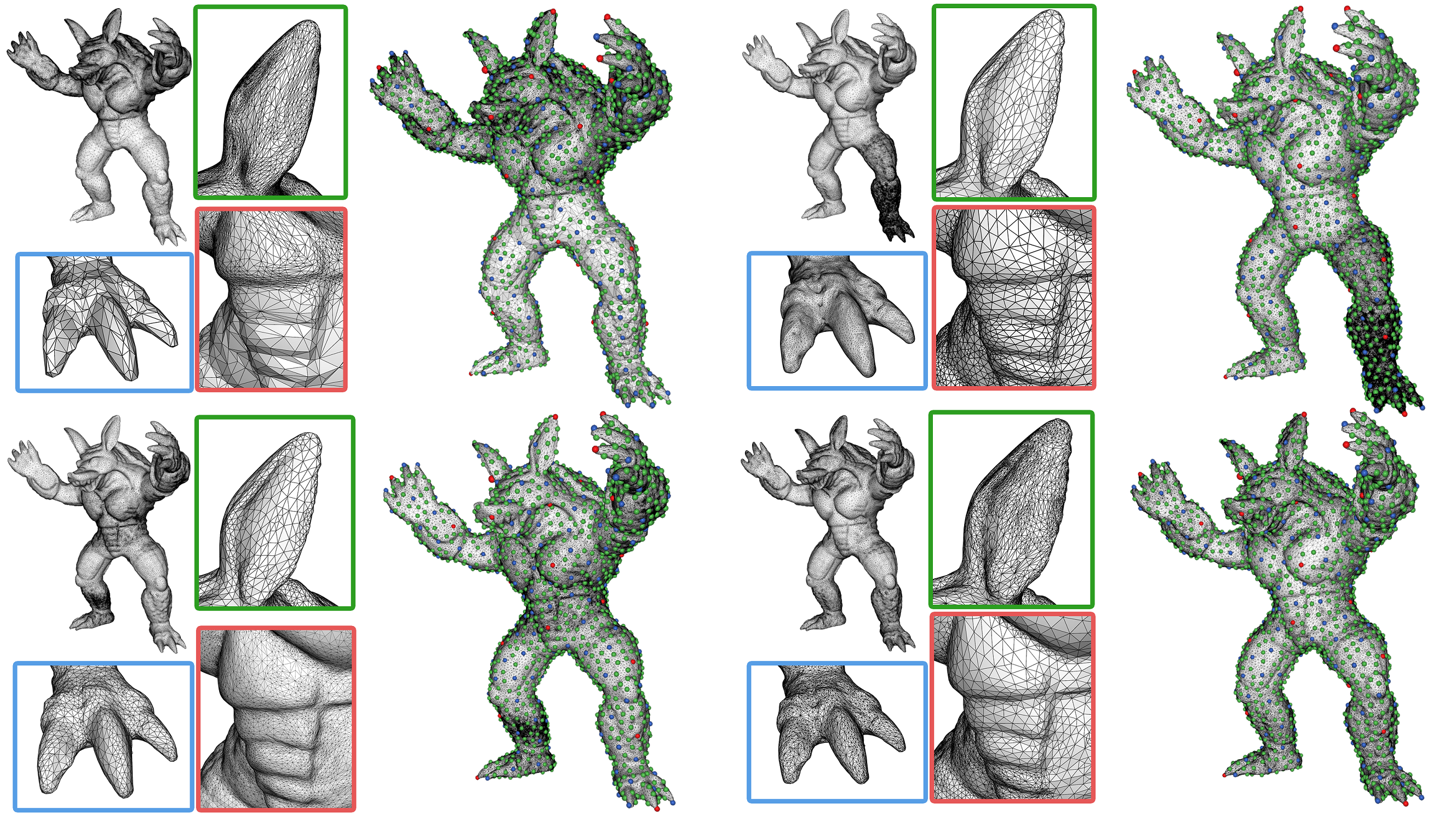}
    \caption{Results of vertex hierarchy construction using farthest point sampling are shown for four meshes that approximate the same surface but have different  spatially varying sampling densities.}
    \label{fig:sampling-adaptive} 
\end{figure}
\subsection{Distance computation}
The construction of the basis functions, see Equation (11) in \cite{Nasikun2021}, requires the computation of geodesic distances. For the evaluation of HSIM, we used  Dijkstra's algorithm on the weighted edge graph of the mesh using the edge lengths as weights. 
Since the basis functions have local support, we stop the single source Dijkstra computation when all vertices in the support of the basis function have been processed. 
Alternatives to Dijkstra's algorithm are the Short Term Vector Dijkstra (STVD) algorithm ~\cite{Campen2013} and the Heat Method~\cite{crane2013geodesics}.  
Table~\ref{tab:basisType} compares timings and iteration counts obtained by using Dijkstra's algorithm, the STVD algorithm and the heat method for basis construction. 
One can see that the required numbers of iterations are similar for all three methods with some slight variations. 
Therefore, the timings for the case that Dijkstra's algorithm is used are comparable to those when the STVD algorithm is used. There is only a small overhead resulting from the additional computational effort of the STVD algorithm compared to Dijkstra's algorithm.  
The heat method is much slower since for each point the distance to all other points is computed instead of only in a local neighborhood.  We would like to note that there are also possibilities to localize the distance computation with the heat method~\cite{Herholz2017}. This, however, would be beyond the scope of this experiment. 
These results illustrate our impression that the STVD algorithm or a localized version of the heat method can be used as alternatives for the basis construction. Since we did not observe any substantial advantages of STVD or the Heat method over Dijkstra's algorithm in our experiments, and to keep the method simple, we used Dijkstra's algorithm for our evaluation of HSIM.

\begin{figure}[]
    \subfloat[HSIM, 10k vertices \label{fig:error-hsim-10k-irreg}]{
        \includegraphics[width=0.49\linewidth]{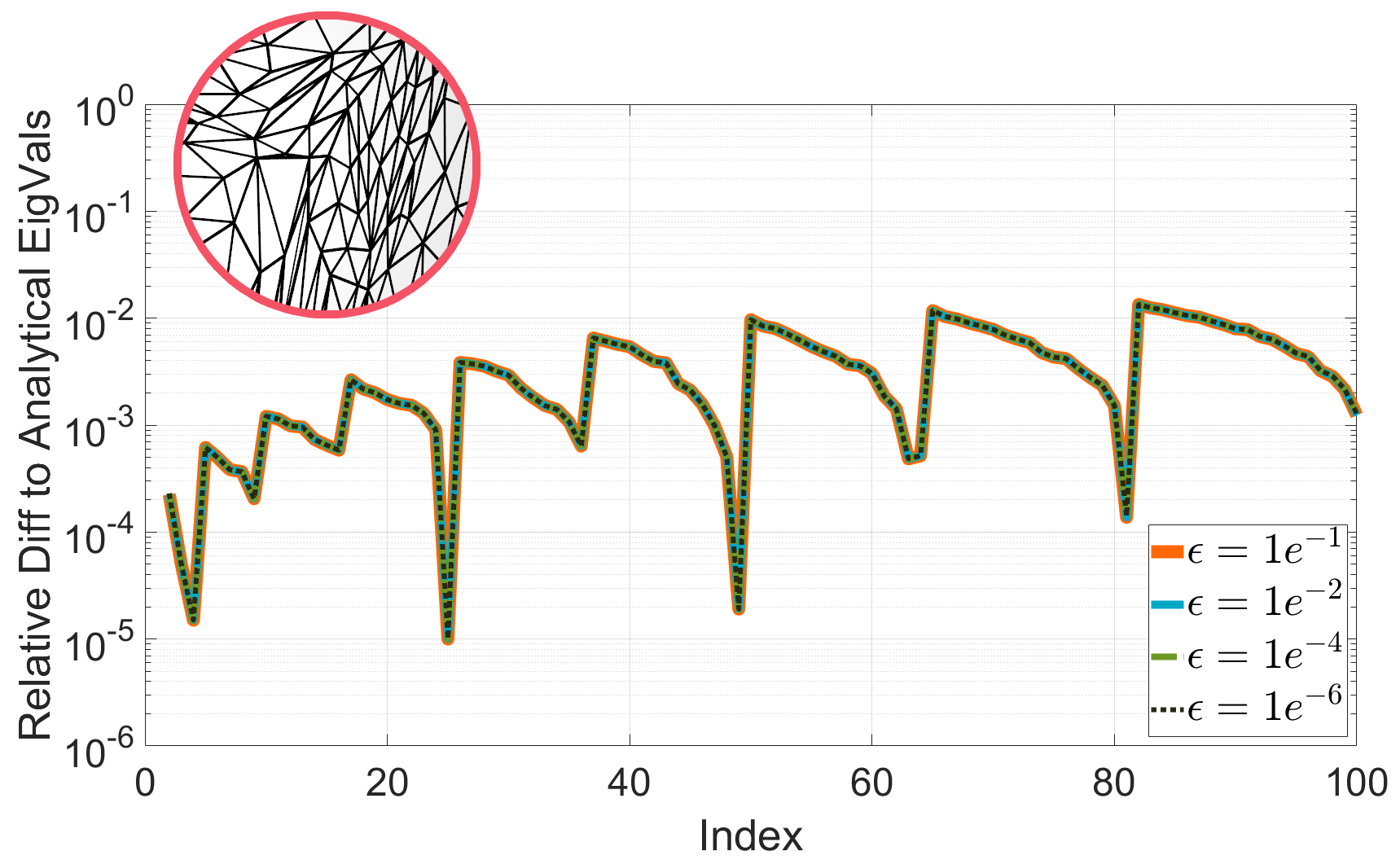}
    }
    \subfloat[HSIM, 100k vertices\label{fig:error-hsim-100k-irreg}]{
        \includegraphics[width=0.47\linewidth]{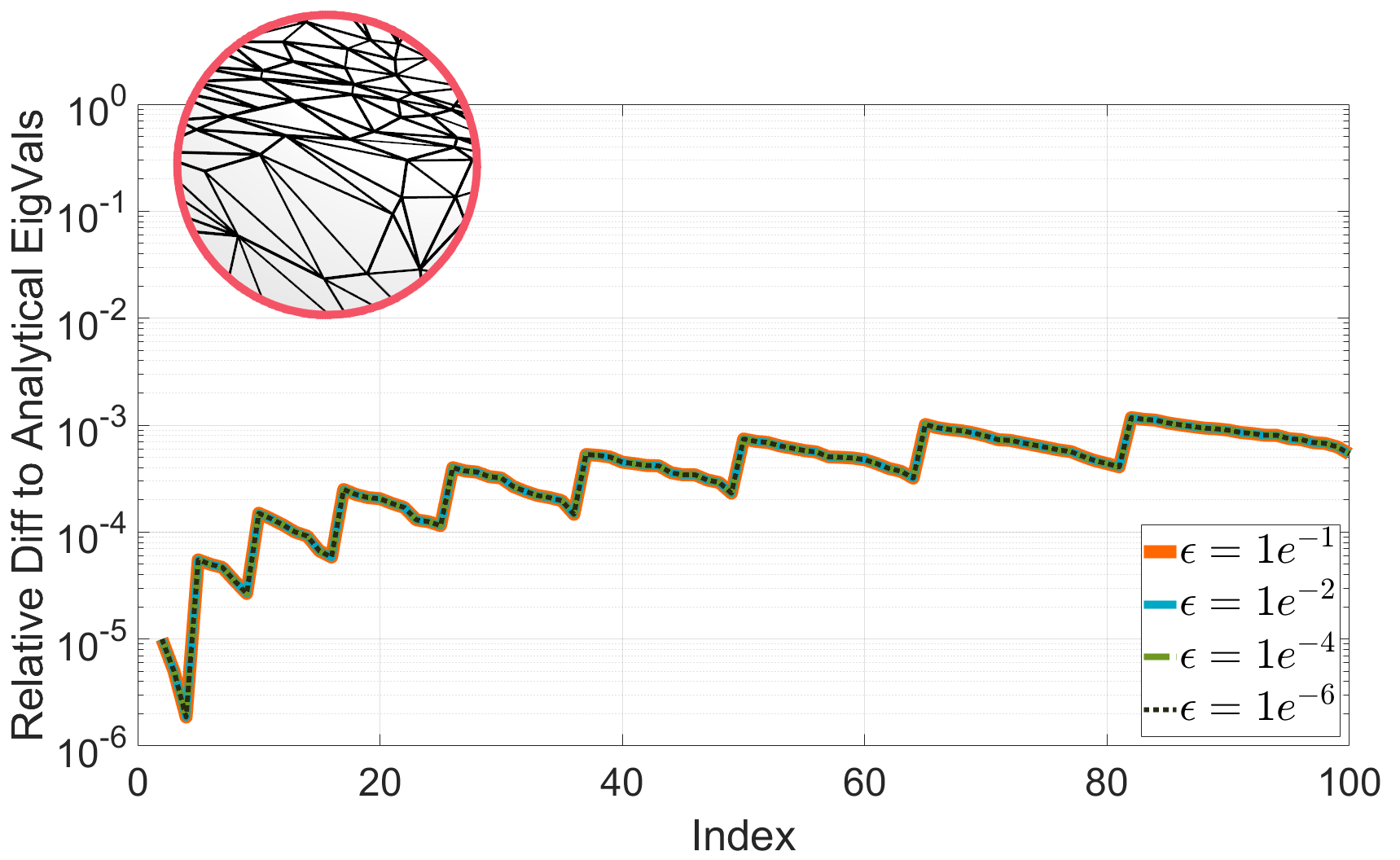}
    }
    \\
    \subfloat[HSIM, 1m vertices\label{fig:error-hsim-1000k-irreg}]{
        \includegraphics[width=0.49\linewidth]{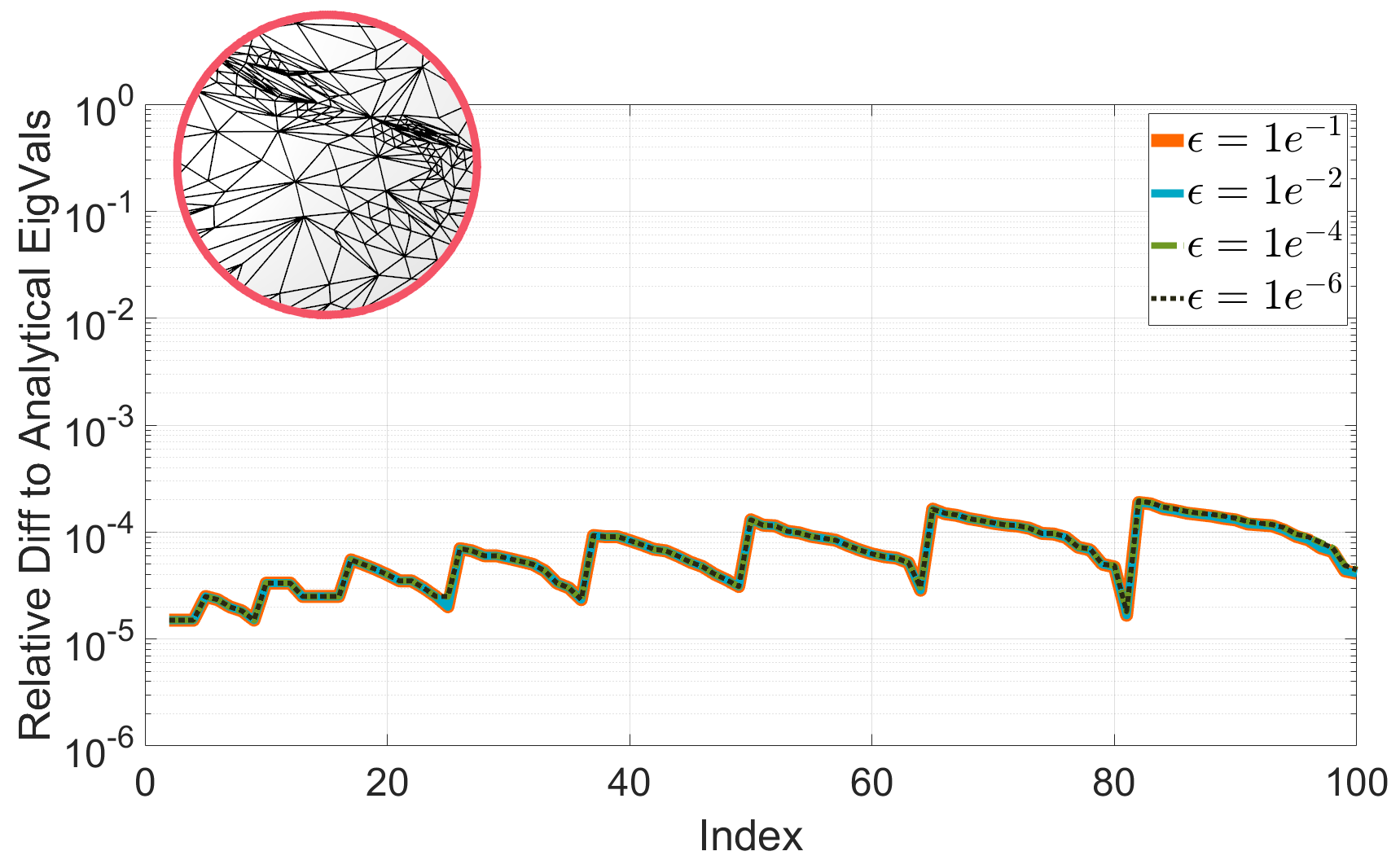}
    }
    \subfloat[IPM, 100k vertices\label{fig:error-ipm-100k-irreg}]{
        \includegraphics[width=0.47\linewidth]{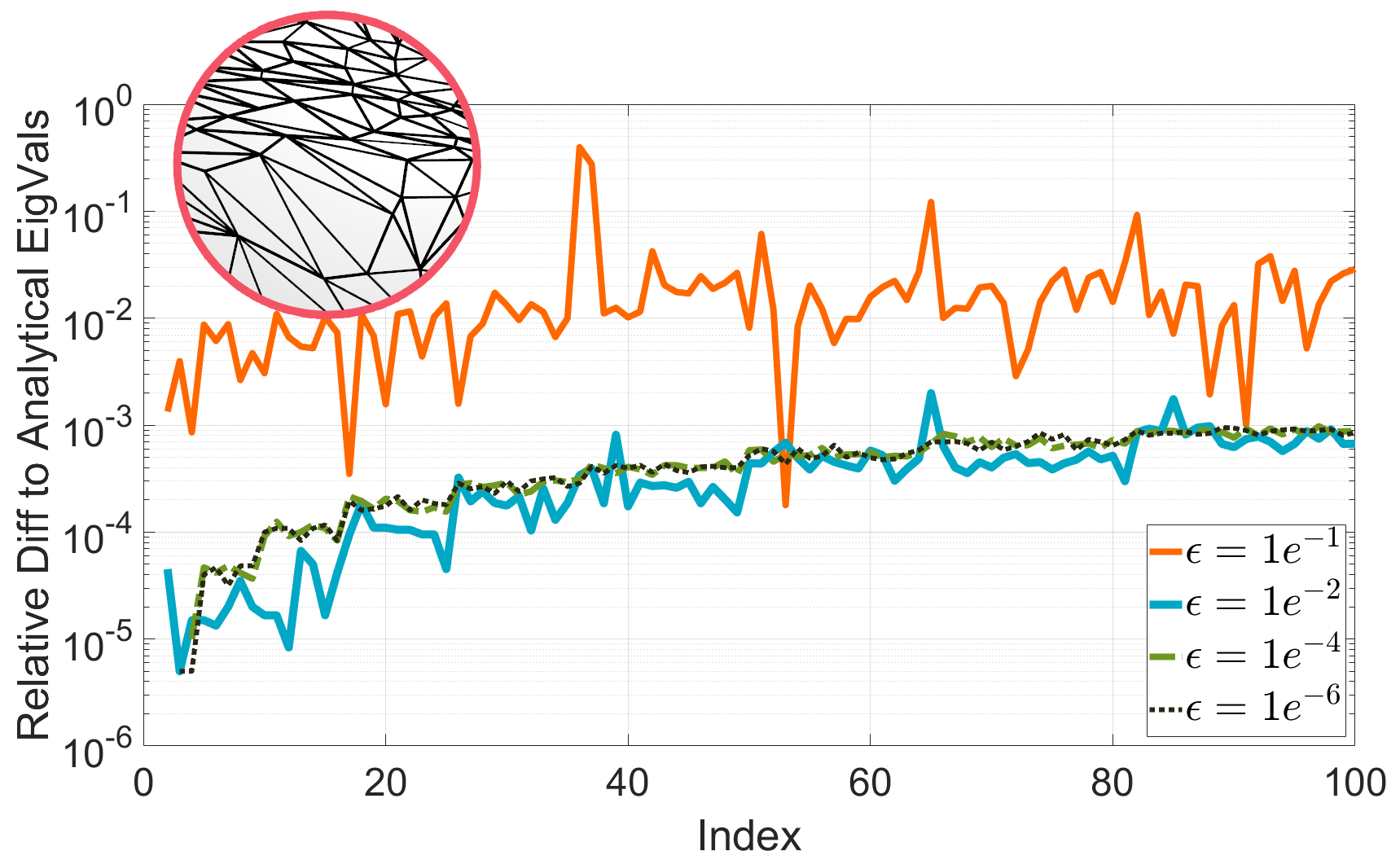}
    }
     \caption{Relative difference of numerical approximations of the eigenvalues
of the unit sphere to the analytic solutions.}
     \label{fig:irreg-sphere} 
\end{figure}

\begin{figure}[]
    \subfloat[Eigenvalues\label{fig:eFlip-TRex-eval}]{
        \includegraphics[width=0.49\linewidth]{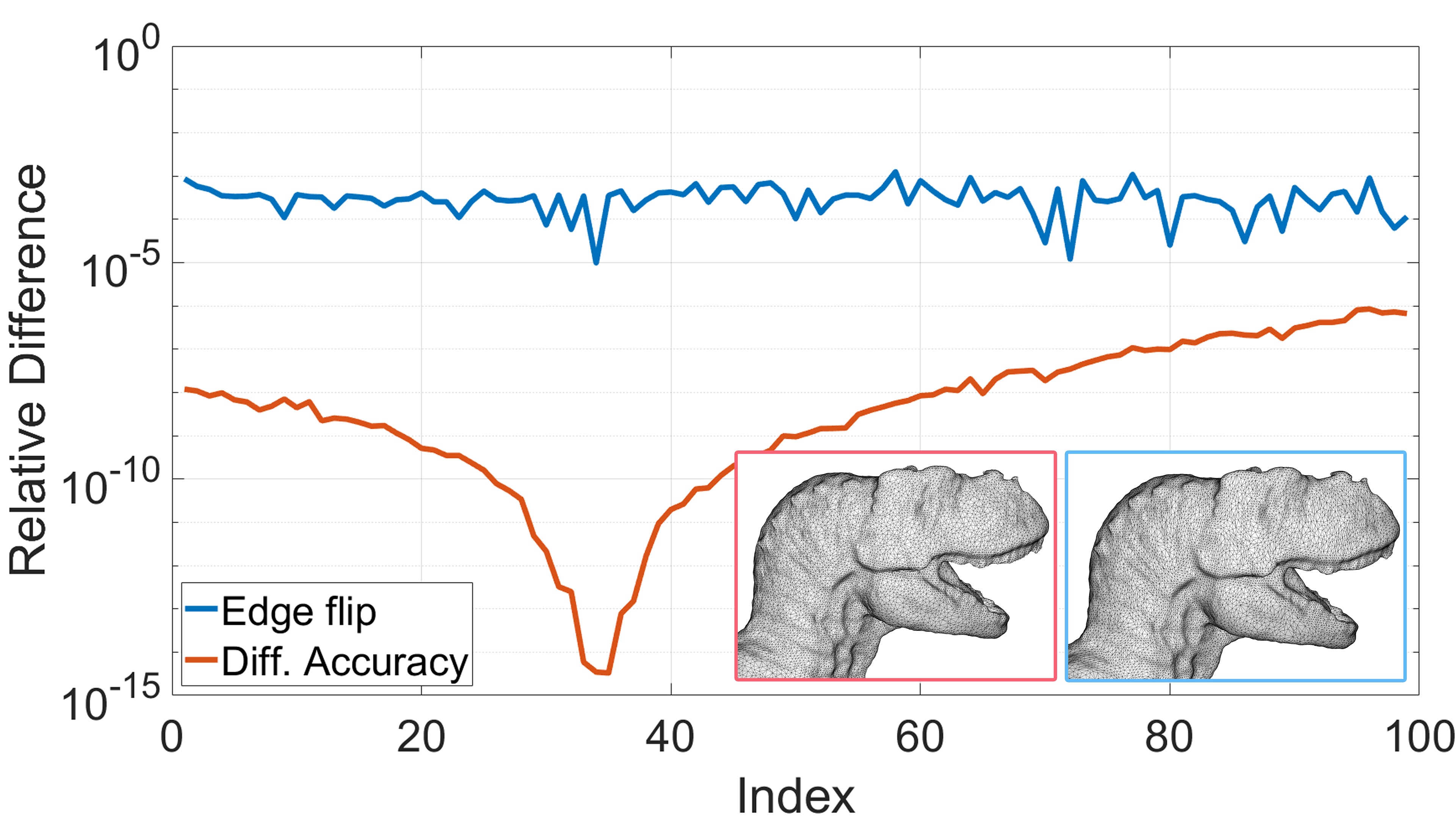}
    }
    \subfloat[Eigenvectors \label{fig:eFlip-TRex-evect}]{
        \includegraphics[width=0.49\linewidth]{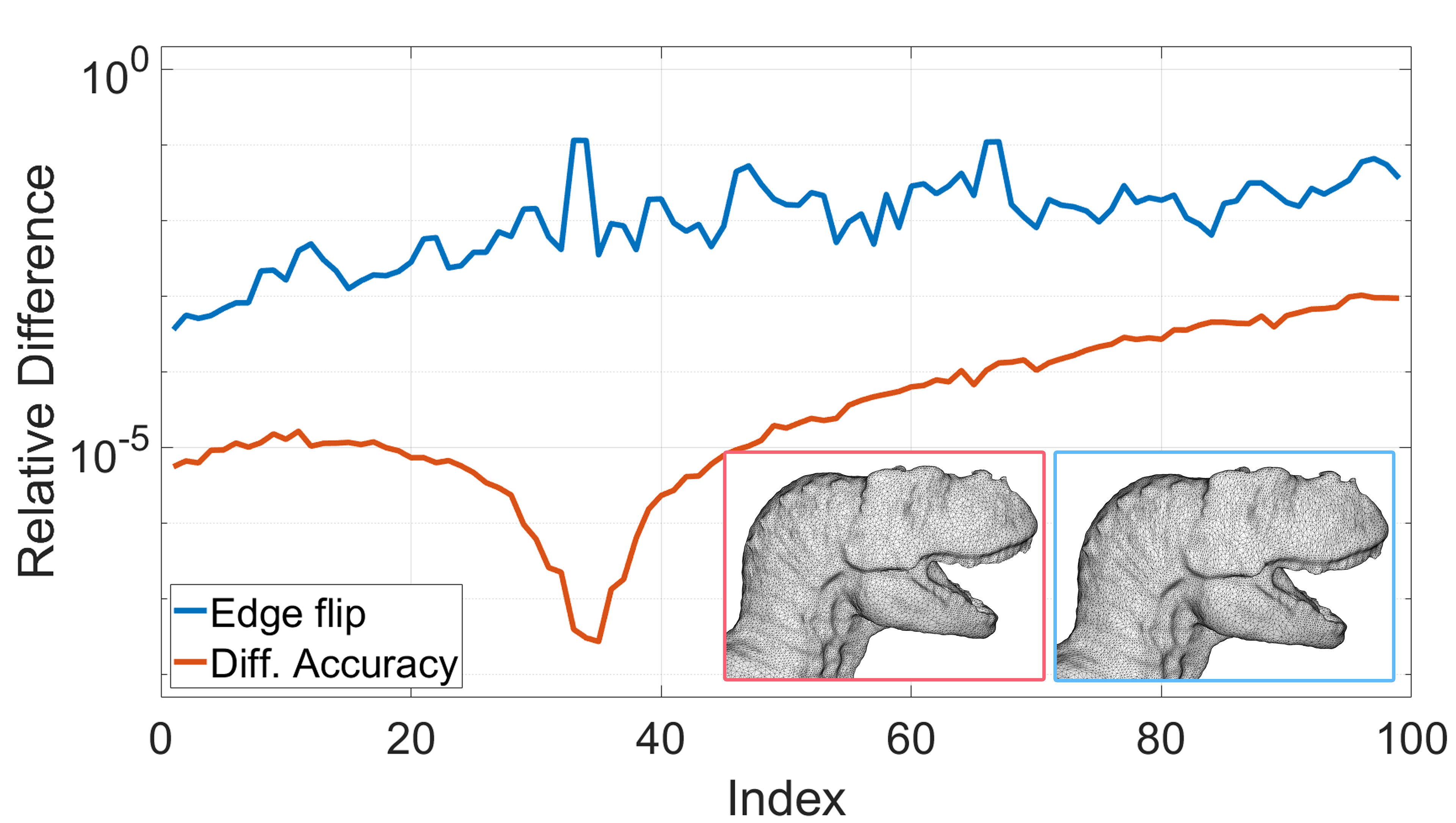}
    }
     \caption{Comparisons of the relative difference of the eigenvalues and the
eigenvectors between two meshes that approximate the same surface (blue
graph) and solutions for different tolerance on one of the meshes (red
graph).}
     \label{fig:edge-flip-TRex} 
\end{figure}

\begin{figure}[]
    \subfloat[Bimba, 100k\label{fig:edge-decimation-1}]{
        \includegraphics[width=0.49\linewidth]{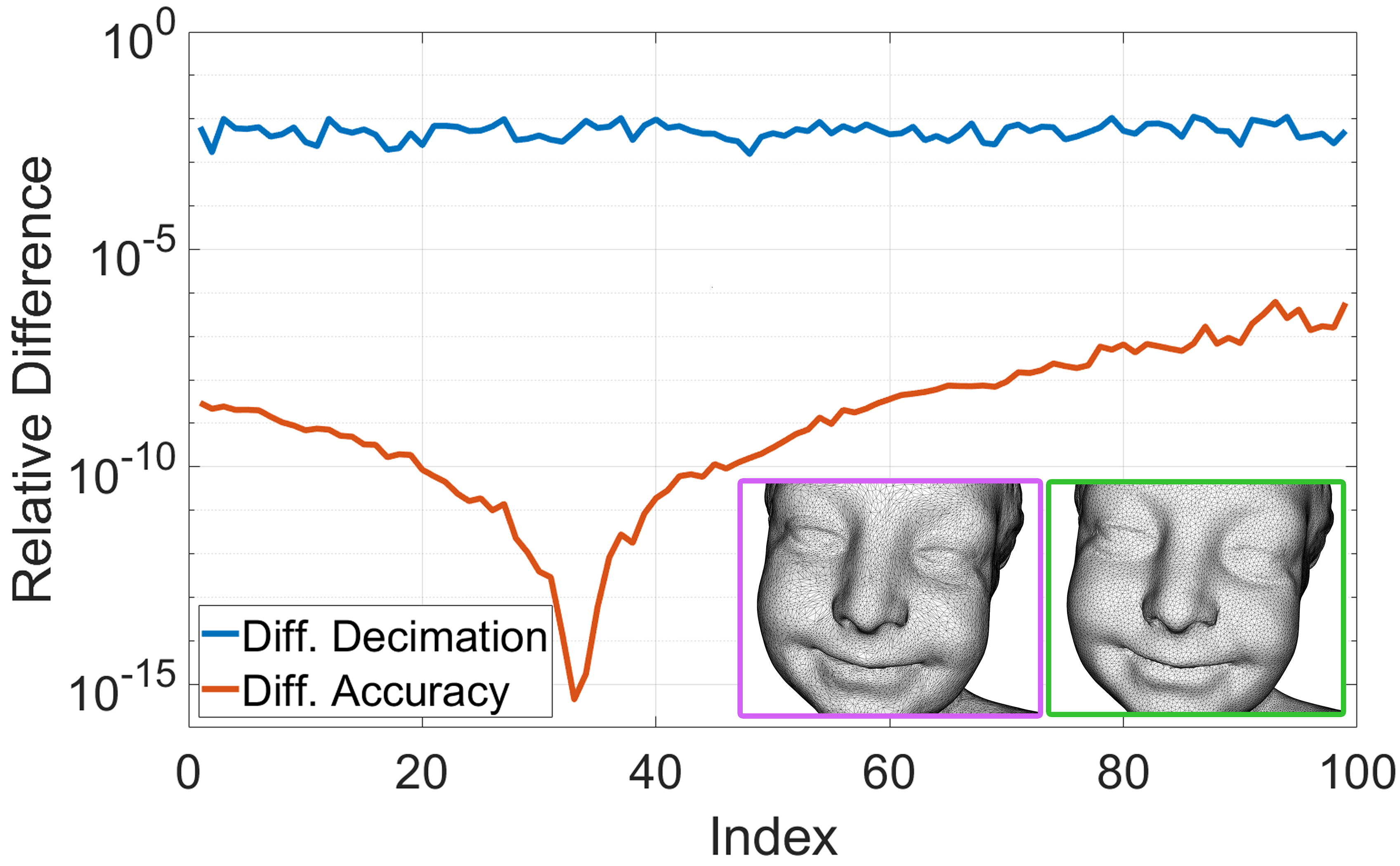}
    }
     \subfloat[Nefertiti 100k\label{fig:decimated-nefertiti-11}]{
        \includegraphics[width=0.47\linewidth]{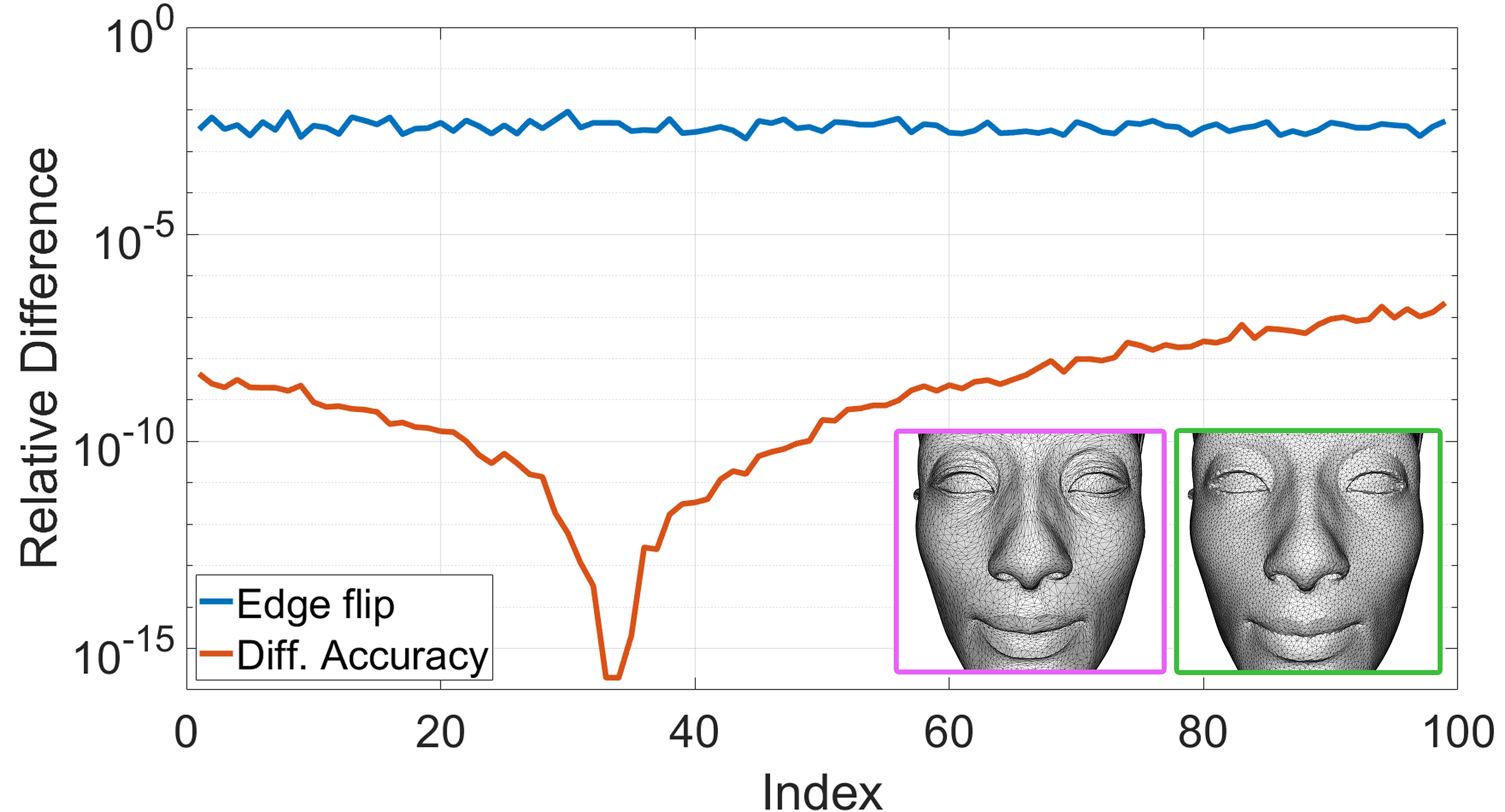}
    }
    \\
        \subfloat[Ramses 100k\label{fig:decimated-ramses-11}]{
        \includegraphics[width=0.49\linewidth]{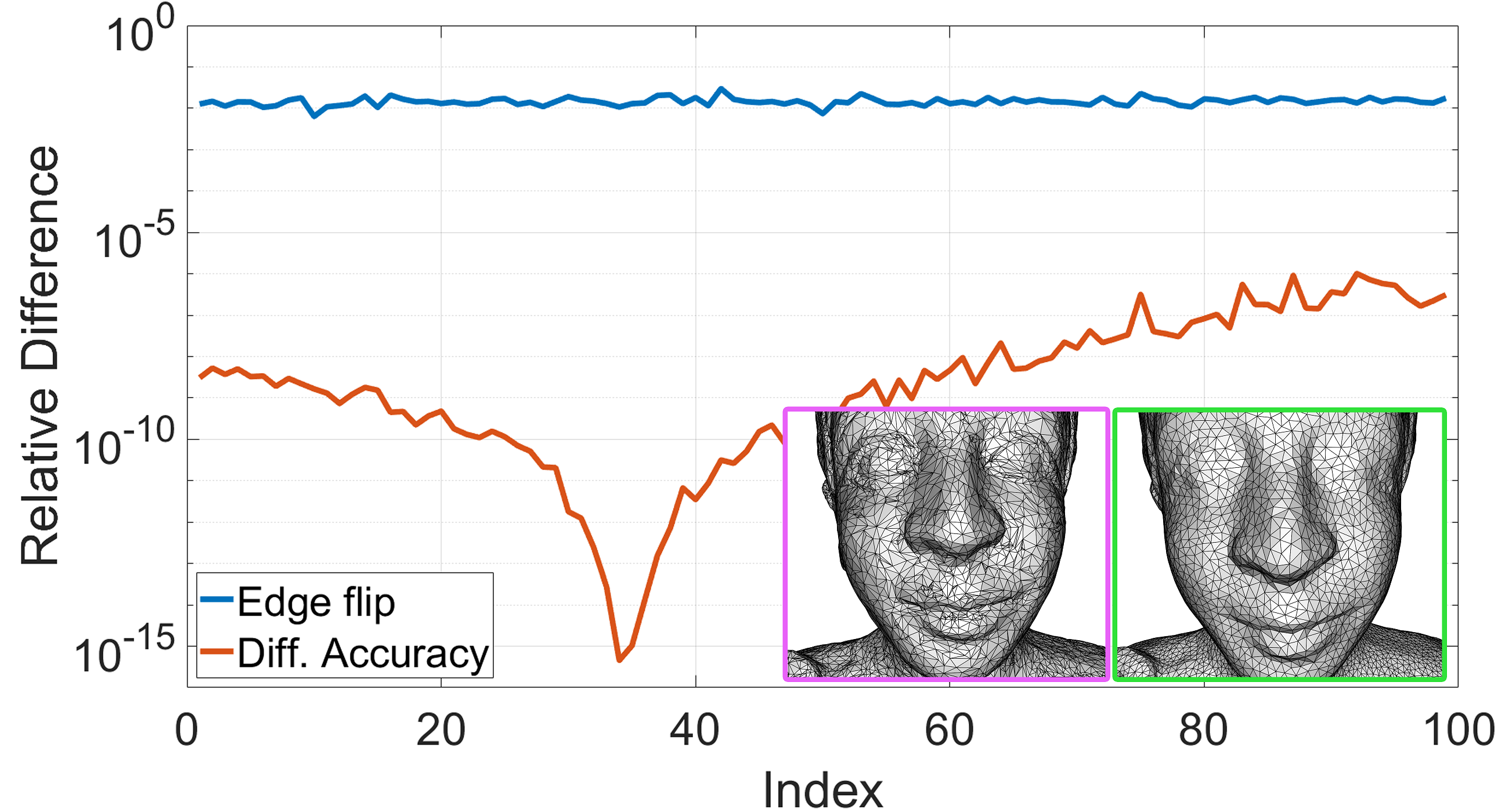}
    }
       \subfloat[Ramses, 500k.\label{fig:edge-decimation-2}]{
        \includegraphics[width=0.49\linewidth]{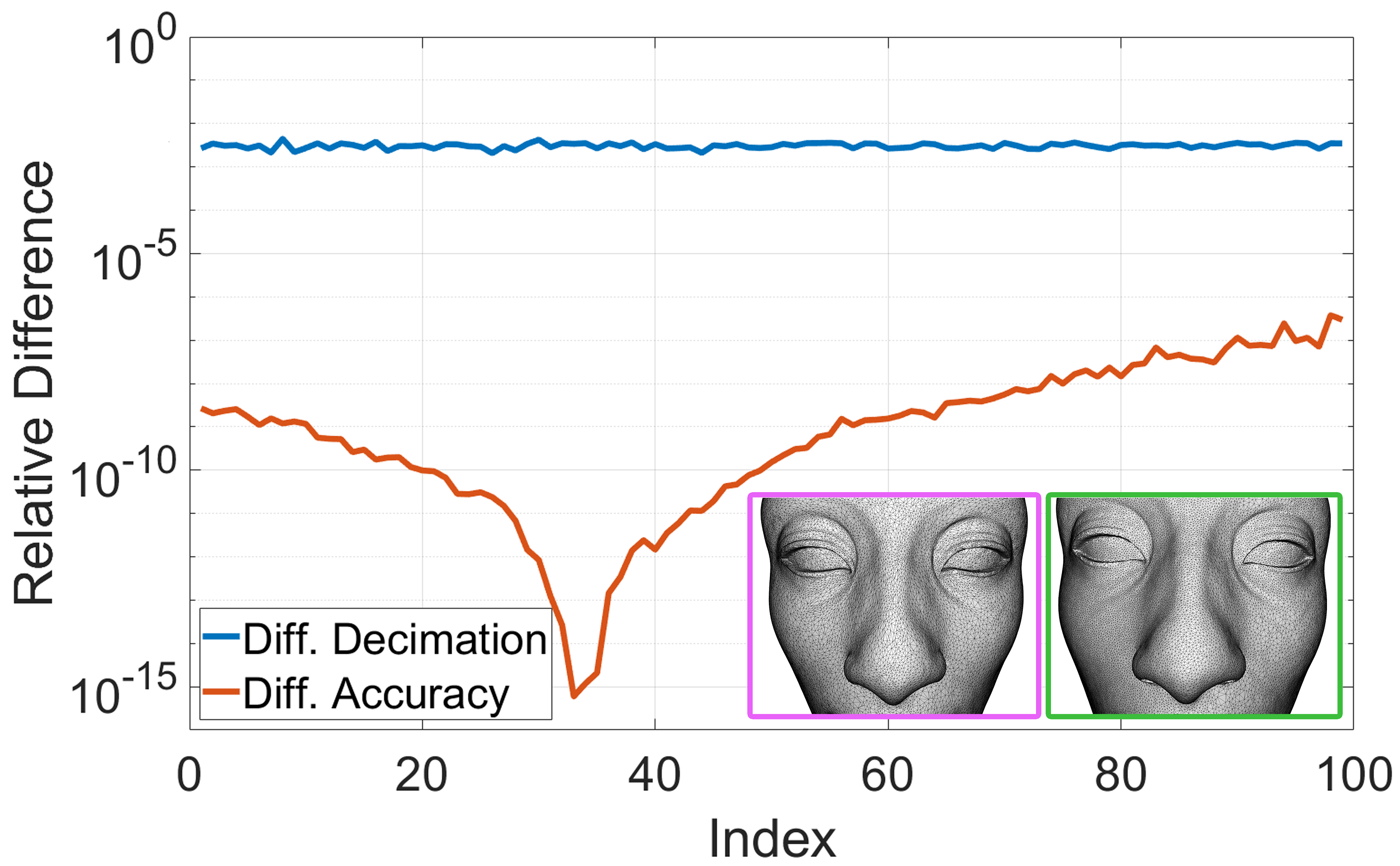}
    }
     \caption{{For four pairs of meshes, each pair approximating the same surface, comparisons of the relative differences between the eigenvalues of the two meshes (blue graphs) and solutions for different convergence tolerance on one of the meshes (red graphs) are shown.}}
     \label{fig:decimation} 
\end{figure}

\subsection{Sampling method}
The nested function spaces we use for HSIM are constructed from a vertex hierarchy, which assigns a level to every vertex of the mesh. 
We use farthest point sampling for computing the distribution of the vertices. In this section, we show examples of vertex hierarchies on different meshes to illustrate why we think farthest point sampling is suitable for this process. 
In the examples, we use meshes with spatially varying resolution, in some areas of the surfaces the triangles are much smaller than in others. 
Figure~\ref{fig:sampling} shows vertex hierarchies on four surfaces. The hierarchies shown have four levels, color-coded in red, blue, and green, with the finest levels encompassing all vertices and not shown. It can be seen that the vertices in the different levels are distributed fairly uniformly over the surface despite the irregular mesh.   
In Figure~\ref{fig:sampling-adaptive} four more examples are shown. In this case, we show vertex hierarchies with the same numbers of vertices in each level on different meshes that approximate the same surface. 
As illustrated in the shown results, we consider the farthest point sampling as a suitable method to build up the vertex hierarchies for HSIM. 

\begin{figure}[t]
    \includegraphics[width=\linewidth]{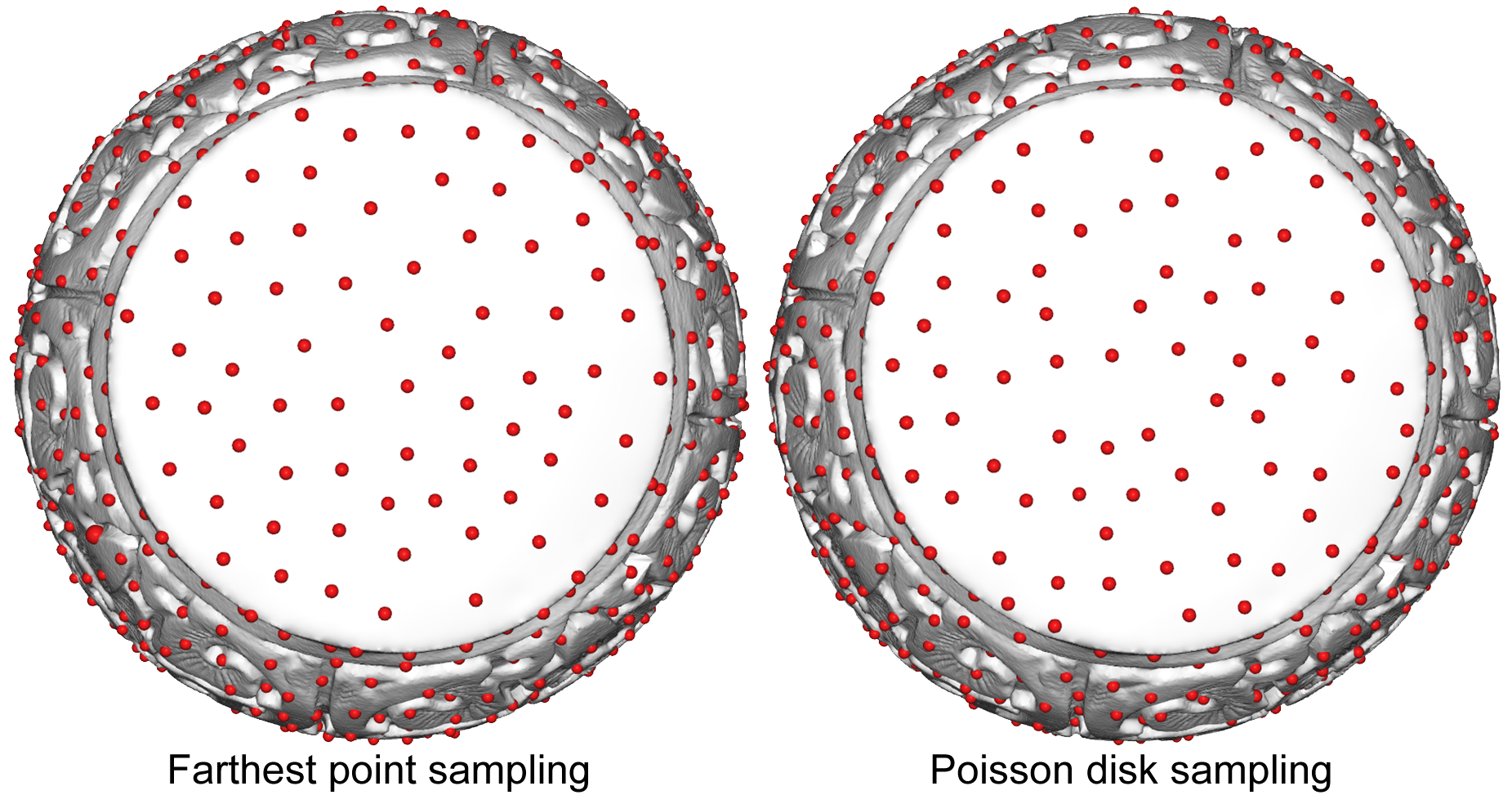}
     \caption{\minorrev{Samplings we computed with farthest point sampling (left) and Poisson disk sampling (right).}}
     \label{fig:FPSvsPDS} 
\end{figure}

\minorrev{
Alternative sampling schemes such as Poisson disk sampling \cite{corsini2012efficient}, which is used in~\cite{Nasikun2018} for the construction of function spaces, could be used instead of farthest point sampling. 
Figure~\ref{fig:FPSvsPDS} shows a comparison of samplings produced with farthest point sampling and Poisson disk sampling. 
A potential benefit of Poisson disk sampling is that the sampling step could be accelerated. Table~\ref{tab:PDS} compares timings required to compute Poisson disk and farthest point samplings of different size on different meshes.
On the other hand, Poisson disk sampling also has drawbacks compared to farthest point sampling. 
While farthest point sampling measures intrinsic distances in the surface, the Poisson disk sampling we compare to uses distances in ambient space, which may deviate when coarser samplings are computed. Another point is that farthest point samplings can be easily split into a hierarchy of levels. For Poisson disk sampling, a division would have to be computed. Lastly, the control parameter for Poisson disk sampling is the disk size rather than the number of samples and our hierarchy construction prescribes the number of samples in each level. We show a visual comparison of the point samplings produced by both sampling methods in Figure~\ref{fig:FPSvsPDS}.
}

\subsection{Convergence tolerance}

This paragraph includes further experiments related to the discussion of the convergence tolerance that we used for the evaluation of HSIM, see also Section 5 of \cite{Nasikun2021}. 
Figure~\ref{fig:irreg-sphere} shows a variant of Figure~\refFigErrorAnalytic~from \cite{Nasikun2021}, where we consider non-regular meshes inscribed to the sphere.  
Figure~\ref{fig:edge-flip-TRex} shows a variant of Figure~\refFigErrFlip~from \cite{Nasikun2021}, using a different mesh. 
 
In Figure~\ref{fig:decimation}, we show results of an additional experiment. As in Figure~\ref{fig:edge-flip-TRex}, we computed eigenpairs on two meshes with tolerances $\varepsilon=10^{-2}$ and as reference with $\varepsilon=10^{-8}$.
We generated the meshes by simplifying one mesh with two different mesh coarsening algorithms. 
We used the Bimba mesh with 500k vertices to get two simplified meshes with 100k vertices each, the Nefertiti mesh with 1m vertices to obtain two simplified meshes with 100k vertices and the Ramses mesh with 750k vertices to get two meshes with 500k vertices and two meshes with 100k vertices. 
In Figure~\ref{fig:decimation}, we plot for all four pairs of meshes the differences between the reference results that are computed with a tolerance of $\varepsilon=10^{-8}$ on both meshes and for one mesh the difference between the results for $\varepsilon=10^{-2}$ and $\varepsilon=10^{-8}$.
For all four pairs of meshes, the difference between the reference results on the two meshes is much larger than the difference between the results for $\varepsilon=10^{-2}$ and $\varepsilon=10^{-8}$.

\begin{table}[t]
\begin{tabular}{|l|r|r|r|r|}
\hline
\multicolumn{1}{|c|}{{  Model (\#Vertices)}}                                 & \multicolumn{1}{c|}{{  \#Eigs}} & \multicolumn{1}{c|}{{  \#Samples}} & \multicolumn{1}{c|}{{  FPS}} & \multicolumn{1}{c|}{{  PDS}} \\ \hline
{  }                                                            & {  50}                         & {  1k}                            & {  1.08}                     & {  0.07}                     \\ \cline{2-5} 
{  }                                                            & {  250}                        & {  11.7k}                         & {  2.01}                     & {  0.18}                     \\ \cline{2-5} 
\multirow{-3}{*}{{  Kitten (137k)}}                             & {  1000}                       & {  16.6k}                         & {  2.37}                     & {  0.24}                     \\ \hline
{  }                                                            & {  50}                         & {  1k}                            & {  1.37}                     & {  0.11}                     \\ \cline{2-5} 
{  }                                                            & {  250}                        & {  14.1k}                         & {  2.93}                     & {  0.24}                     \\ \cline{2-5} 
\multirow{-3}{*}{{  Vase-Lion (200k)}}                          & {  1000}                       & {  200k}                          & {  3.42}                     & {  0.32}                     \\ \hline
{  }                                                            & {  50}                         & {  1k}                            & {  4.25}                     & {  0.27}                     \\ \cline{2-5} 
{  }                                                            & {  250}                        & {  21.9k}                         & {  9.96}                     & {  0.54}                     \\ \cline{2-5} 
\multirow{-3}{*}{{  Knot-Stars (450k)}}                         & {  1000}                       & {  30.9k}                         & {  11.98}                    & {  0.68}                     \\ \hline
{  }                                                            & {  50}                         & {  1k}                            & {  5.05}                     & {  0.38}                     \\ \cline{2-5} 
{  }                                                            & {  250}                        & {  23.9k}                         & {  12.12}                    & {  0.64}                     \\ \cline{2-5} 
\multirow{-3}{*}{{  Oilpump (570k)}}                            & {  1000}                       & {  33.8k}                         & {  14.59}                    & {  0.77}                     \\ \hline
\multicolumn{1}{|r|}{{  }}                                      & {  50}                         & {  1k}                            & {  5.09}                     & {  0.46}                     \\ \cline{2-5} 
\multicolumn{1}{|r|}{{  }}                                      & {  250}                        & {  26.5k}                         & {  14.69}                    & {  0.71}                     \\ \cline{2-5} 
\multicolumn{1}{|r|}{\multirow{-3}{*}{{  Red-Circular (700k)}}} & {  1000}                       & {  37.5k}                         & {  17.64}                    & {  0.91}                     \\ \hline
\end{tabular}
\caption{\minorrev{Comparison of computation time (in seconds) of farthest point sampling (FPS) and Poisson-disk sampling (PDS).}}
	\label{tab:PDS}
\end{table}

\begin{table}[b]
\centering
\resizebox{0.85\columnwidth}{!}{
\begin{tabular}{|l|r|r|r|r|r|}
\hline
\multicolumn{1}{|c|}{Model}     & \multicolumn{1}{c|}{\#Verts.} & \#Eigs. & HSIM   & \textsc{Matlab} & LOBPCG     \\ \hline
\multirow{3}{*}{Gargoyle}       & \multirow{3}{*}{85k}          & 50      & 4.2    & 5.1    & 27.1       \\ \cline{3-6} 
                                &                               & 250     & 19.9   & 49.3   & 154.7      \\ \cline{3-6} 
                                &                               & 1000    & 88.6   & 442.3  & 995.0      \\ \hline
\multirow{3}{*}{Chinese Dragon} & \multirow{3}{*}{135k}         & 50      & 7.5    & 7.6    & 57.7       \\ \cline{3-6} 
                                &                               & 250     & 30.5   & 56.3   & 298.0      \\ \cline{3-6} 
                                &                               & 1000    & 127.5  & 523.3  & 1678.7     \\ \hline
\multirow{3}{*}{Dragon}         & \multirow{3}{*}{150k}         & 50      & 7.2    & 9.6    & 83.0       \\ \cline{3-6} 
                                &                               & 250     & 36.5   & 65.7   & 325.1      \\ \cline{3-6} 
                                &                               & 1000    & 143.8  & 795.9  & 2102.5     \\ \hline
\multirow{3}{*}{Blade}          & \multirow{3}{*}{200k}         & 50      & 10.5   & 14.1   & 97.9       \\ \cline{3-6} 
                                &                               & 250     & 49.2   & 93.9   & 453.6      \\ \cline{3-6} 
                                &                               & 1000    & 177.8  & 1091.9 & 2591.0     \\ \hline
\multirow{3}{*}{Fertility}      & \multirow{3}{*}{240k}         & 50      & 14.6   & 16.6   & 133.3      \\ \cline{3-6} 
                                &                               & 250     & 90.8   & 121.6  & 678.6      \\ \cline{3-6} 
                                &                               & 1000    & 236.1  & 1369.9 & 4003.0     \\ \hline
\multirow{3}{*}{Rocker-Arm}     & \multirow{3}{*}{270k}         & 50      & 17.5   & 22.2   & 135.9      \\ \cline{3-6} 
                                &                               & 250     & 73.9   & 175.8  & 744.4      \\ \cline{3-6} 
                                &                               & 1000    & 252.1  & 1537.2 & 4837.9     \\ \hline
\multirow{3}{*}{Pulley}         & \multirow{3}{*}{300k}         & 50      & 19.5   & 21.1   & 228.3      \\ \cline{3-6} 
                                &                               & 250     & 85.3   & 222.1  & 837.3      \\ \cline{3-6} 
                                &                               & 1000    & 339.5  & 1795.0 & 5416.9     \\ \hline
\multirow{3}{*}{Eros}           & \multirow{3}{*}{400k}         & 50      & 30.7   & 43.6   & 194.5      \\ \cline{3-6} 
                                &                               & 250     & 127.2  & 267.7  & 1305.5     \\ \cline{3-6} 
                                &                               & 1000    & 322.6  & 2748.4 & Mem. Bound \\ \hline
\multirow{3}{*}{Bimba}          & \multirow{3}{*}{500k}         & 50      & 29.4   & 31.4   & 236.1      \\ \cline{3-6} 
                                &                               & 250     & 132.8  & 254.9  & 1255.3     \\ \cline{3-6} 
                                &                               & 1000    & 569.3  & 3208.0 & Mem. Bound \\ \hline
\multirow{3}{*}{Oilpump}        & \multirow{3}{*}{570k}         & 50      & 41.1   & 46.1   & 310.3      \\ \cline{3-6} 
                                &                               & 250     & 154.2  & 315.6  & 1864.4     \\ \cline{3-6} 
                                &                               & 1000    & 690.9  & 3354.7 & Mem. Bound \\ \hline
\multirow{3}{*}{Rolling stage}  & \multirow{3}{*}{680k}         & 50      & 54.6   & 57.8   & 326.3      \\ \cline{3-6} 
                                &                               & 250     & 197.6  & 386.5  & 2301.2     \\ \cline{3-6} 
                                &                               & 1000    & 891.6  & 4064.1 & Mem. Bound \\ \hline
\multirow{3}{*}{Ramses}         & \multirow{3}{*}{825k}         & 50      & 49.1   & 62.6   & 458.9      \\ \cline{3-6} 
                                &                               & 250     & 221.4  & 413.4  & 2339.9     \\ \cline{3-6} 
                                &                               & 1000    & 1149.1 & 4979.2 & Mem. Bound \\ \hline
\multirow{3}{*}{Nefertiti}      & \multirow{3}{*}{1m}           & 50      & 64.0   & 62.8   & 396.2           \\ \cline{3-6} 
                                &                               & 250     & 305.4  & 682.7  & 2482.2           \\ \cline{3-6} 
                                &                               & 500     & 277.9  & 1654.5 & Mem. Bound           \\ \hline
\end{tabular}
}
\caption{Comparisons of timings of HSIM, \textsc{Matlab}'s Lanczos solver and LOBPCG for Laplace--Beltrami eigenproblems on different meshes. Renderings of the meshes are show in Figure~\ref{fig:models}.}
	\label{tab:comparison}
\end{table}

\begin{figure}[]
    \includegraphics[width=\linewidth]{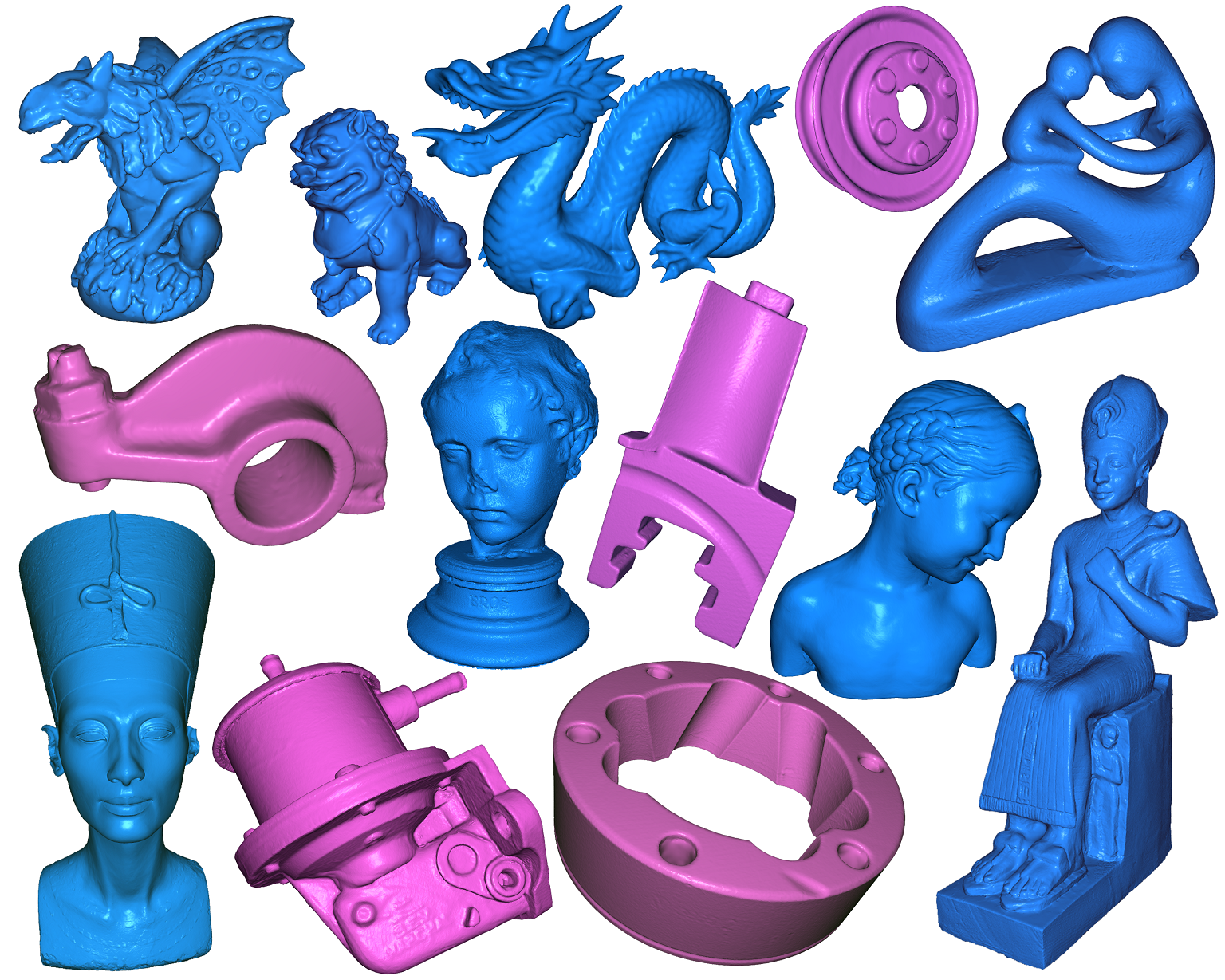}
     \caption{Renderings of the meshes used for the comparisons listed in Table~\ref{tab:comparison}.}
     \label{fig:models} 
\end{figure}

\begin{figure}[b]

  \centering
  \subfloat[32 eigenpairs]{\includegraphics[width=0.47\linewidth]{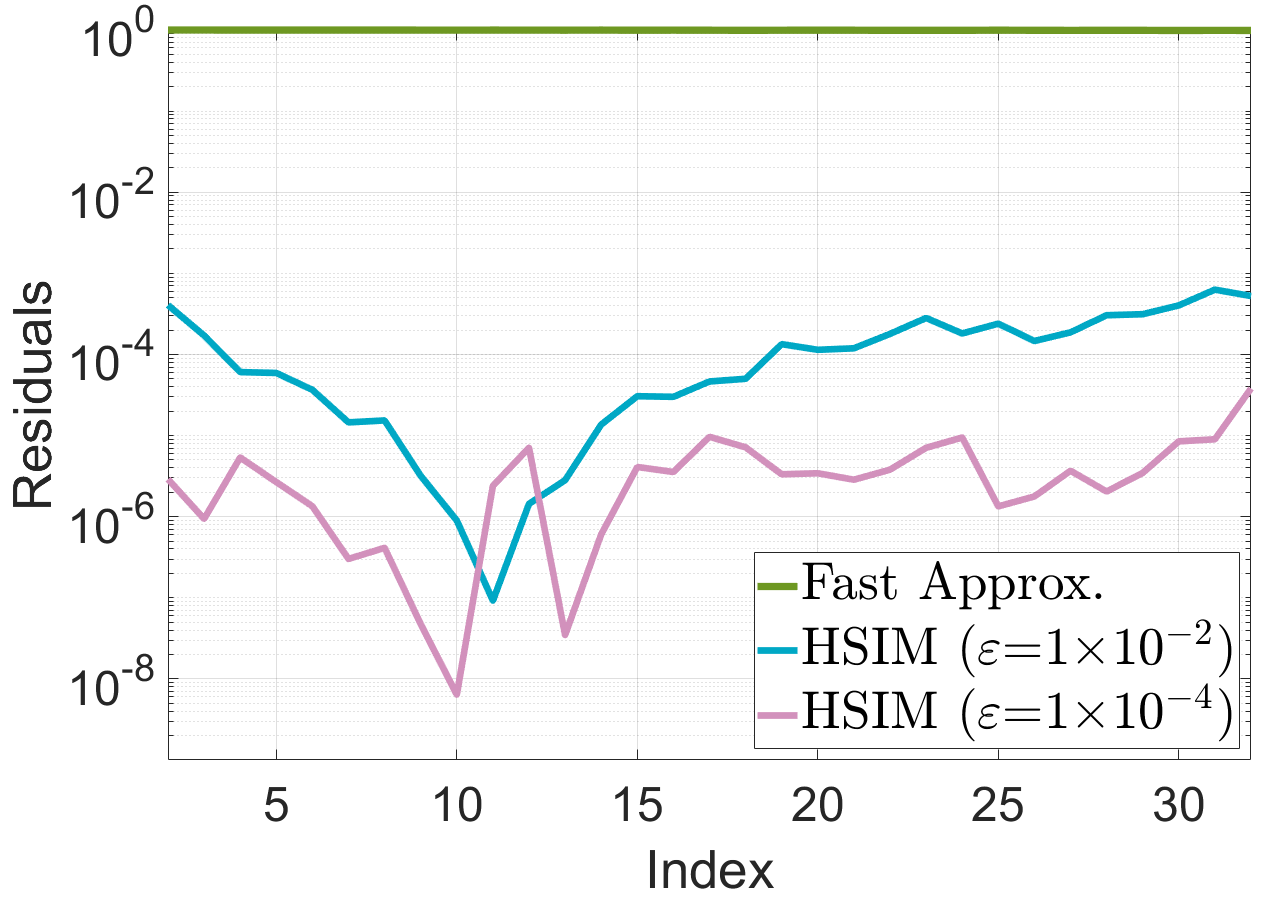}\label{fig:vsSGP11}}
  \hfill
  \subfloat[250 eigenpairs]{\includegraphics[width=0.49\linewidth]{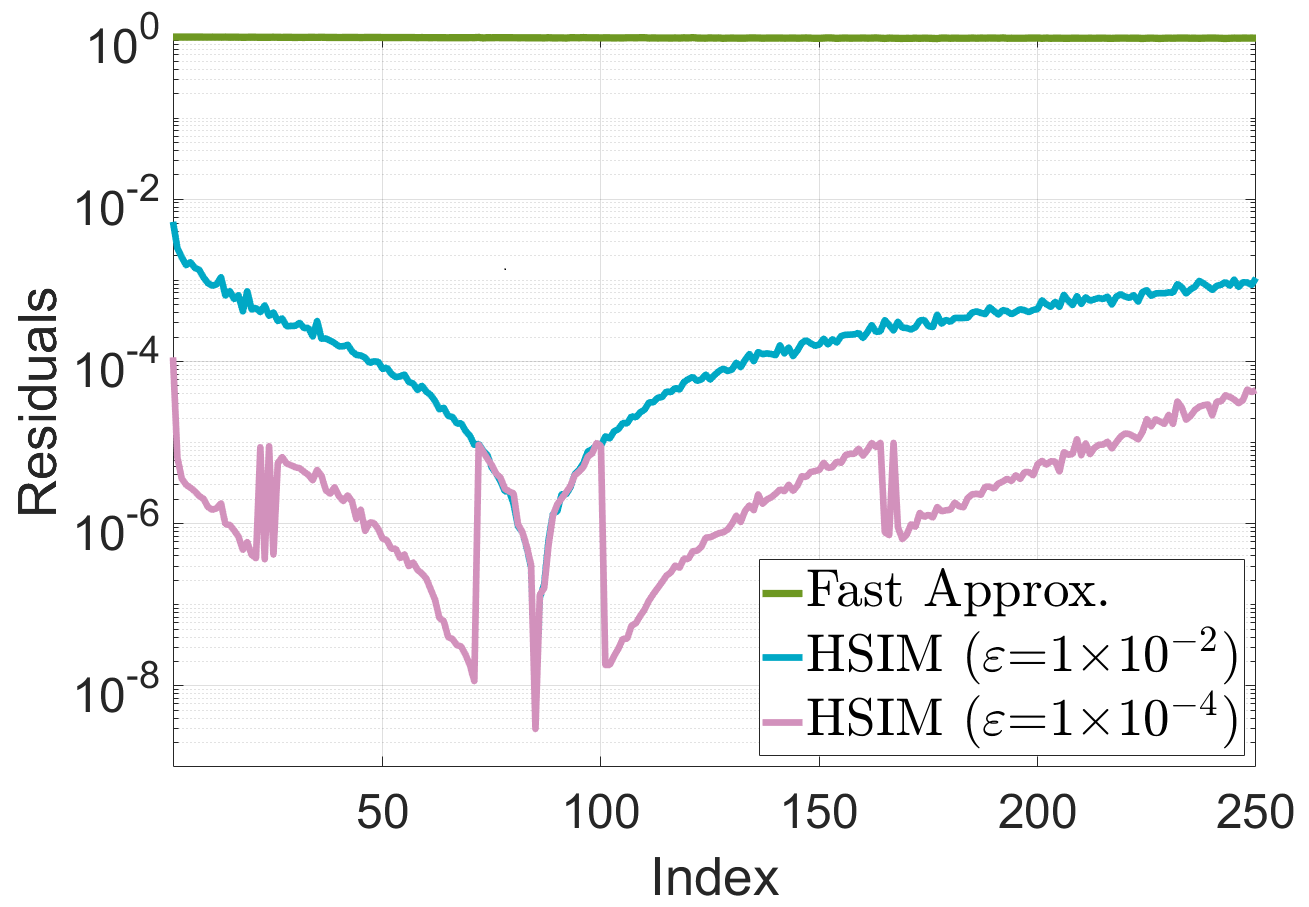}\label{fig:vsSGP12}}  
  \parbox[t]{1.0\columnwidth}{\relax  }
 \\
   \subfloat[32 eigenpairs]{\includegraphics[width=0.47\linewidth]{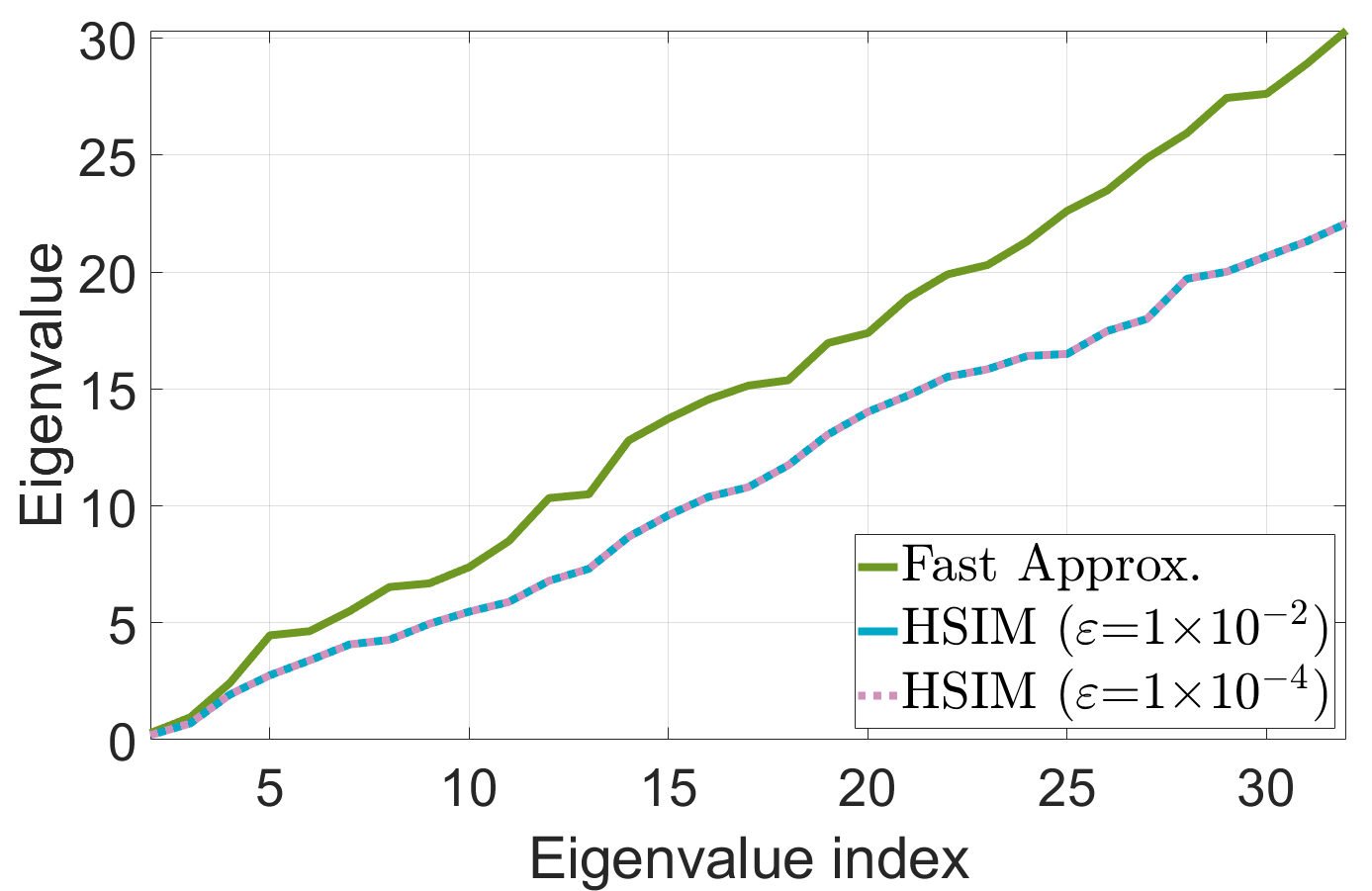}\label{fig:vsSGP21}}
  \hfill
  \subfloat[250 eigenpairs]{\includegraphics[width=0.49\linewidth]{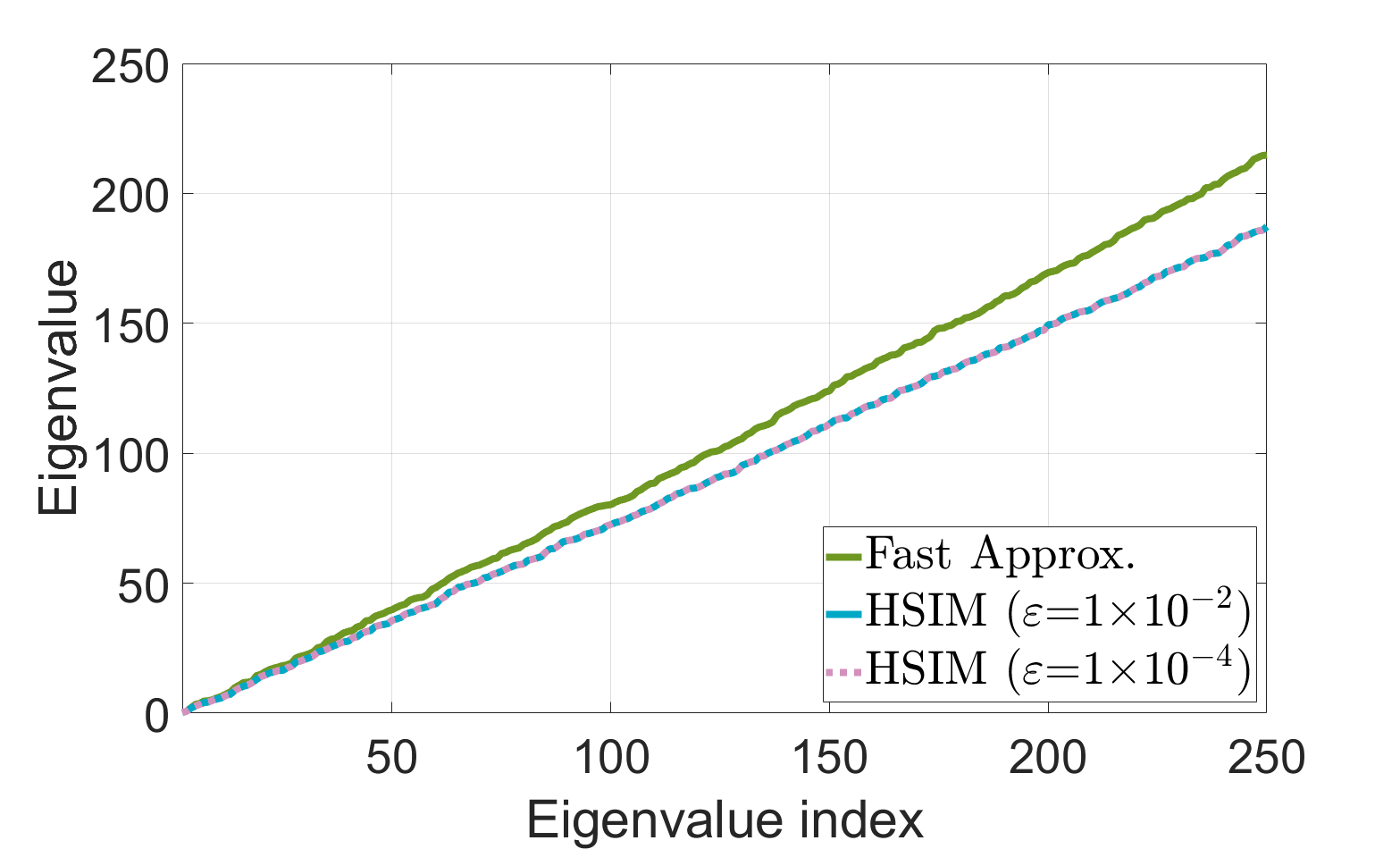}\label{fig:vsSGP22}}  
  \parbox[t]{1.0\columnwidth}{\relax  }
  \caption{
  \label{fig:vsSGP1}
   Top row: Plot of the residuals for the computation of the lowest 32 and 250 Laplace--Beltrami eigenvalues of the Dragon model with 150k vertices. Results for the fast approximation scheme from \cite{Nasikun2018} and HSIM with tolerance $\varepsilon=10^{-2}$ and $\varepsilon=10^{-4}$ are shown.
   Bottom row: The computed eigenvalues are plotted. 
}
\end{figure}

\begin{table*}[tb]
{
\centering
\begin{tabular}{|l|r||r|r|r|r||r|r|r|r||r|r|r|r|}
\hline
\multicolumn{1}{|c|}{\multirow{2}{*}{Model}} & \multicolumn{1}{c||}{\multirow{2}{*}{\#Eigs}} & \multicolumn{4}{|c||}{Laplacian}                                                                                     & \multicolumn{4}{c||}{Hamiltonian (t=0.1)}                                                                                       & \multicolumn{4}{c|}{Hamiltonian (t=1.0)}                                                                                   \\ \cline{3-14} 
\multicolumn{1}{|c|}{}                       & \multicolumn{1}{c|}{}                        & \multicolumn{1}{|c|}{Hier.} & \multicolumn{1}{c|}{Solve} & \multicolumn{1}{c|}{\#Iter} & \multicolumn{1}{c||}{Total} & \multicolumn{1}{c|}{Hier.} & \multicolumn{1}{c|}{Solve} & \multicolumn{1}{c|}{\#Iter} & \multicolumn{1}{c||}{Total} & \multicolumn{1}{c|}{Hier.} & \multicolumn{1}{c|}{Solve} & \multicolumn{1}{c|}{\#Iter} & \multicolumn{1}{c|}{Total} \\ \hline \hline 
\multirow{3}{*}{Cube (25k)}                  & 50                                           & 0.3                        & 2.1                        & F|1                         & 2.4                        & 0.3                        & 2.1                        & F|1                         & 2.4                        & 0.3                        & 2.1                        & F|1                         & 2.4                        \\ \cline{2-14} 
                                             & 250                                          & 0.6                        & 6.7                        & F|1|1                       & 7.3                        & 0.7                        & 6.6                        & F|1|1                       & 7.3                        & 0.7                        & 9.2                        & F|3|1                       & 9.9                        \\ \cline{2-14} 
                                             & 1000                                         & 0.7                        & 67.1                       & F|5|2                       & 66.8                       & 0.8                        & 65.9                       & F|5|2                       & 66.8                       & 0.6                        & 67.3                       & F|5|2                       & 67.9                       \\ \hline
\multirow{3}{*}{Blade (200k)}                & 50                                           & 2.8                        & 7.4                        & F|1                         & 10.2                       & 2.9                        & 7.4                        & F|1                         & 10.3                       & 2.8                        & 10.5                       & F|2                         & 13.3                       \\ \cline{2-14} 
                                             & 250                                          & 7.3                        & 42.8                       & F|2|1                       & 50.1                       & 7.1                        & 42.7                       & F|2|1                       & 49.8                       & 7.2                        & 45.8                       & F|3|1                       & 53.0                       \\ \cline{2-14} 
                                             & 1000                                         & 8.8                        & 158.7                      & F|2|1                       & 167.5                      & 8.8                        & 167.5                      & F|2|1                       & 176.3                      & 8.9                        & 204.0                      & F|4|1                       & 212.9                      \\ \hline
\multirow{3}{*}{Bimba (500k)}               & 50                                           & 7.9                       & 22.6                       & F|1                         & 30.5                       & 8.0                            & 23.8                            & F|1                             & 31.7                            & 7.7                            & 30.9                            & F|2                             & 38.7                            \\ \cline{2-14} 
                                             & 250                                          & 26.4                       & 105.1                      & F|2|1                       & 131.6                      & 26.4                       & 106.8                      & F|2|1                       & 133.2                      & 26.3                       & 121.3                      & F|3|1                       & 147.6                      \\ \cline{2-14} 
                                             & 1000                                         & 34.1                           &  519.5                          & F|3|1 &    553.6                         & 33.5                           &    510.0                                                    & F|3|1                            & 543.5                           & 34.9                           & 643.3                            & F|6|1                            & 678.2                           \\ \hline
\end{tabular}
\caption{Timings and iteration counts for Laplace--Beltrami and Hamiltonian eigenproblems are shown.
	}
	\label{tab:hamiltonian}
}
\end{table*}

\section{Comparisons}
In this section, we show additional comparisons to alternative approaches for solving eigenproblems. 


\begin{figure*}[]
    \includegraphics[width=0.95\linewidth]{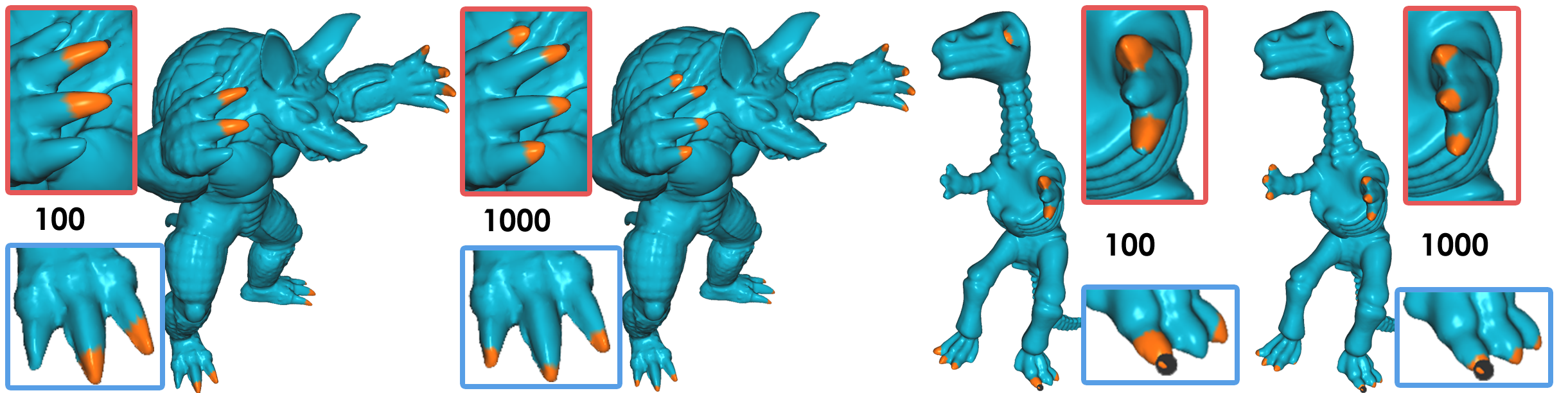}
     \caption{Points similar to a fingertip of the armadillo mesh and to a toe tip for the dinosaur mesh are indicated by binary color-coding. The similarity is computed using the heat kernel distance. Results for heat kernel distance estimation using 100 and 1000 eigenpairs are shown.}
     \label{fig:HKS-region} 
\end{figure*}

\begin{figure*}[]
    \includegraphics[width=0.95\linewidth]{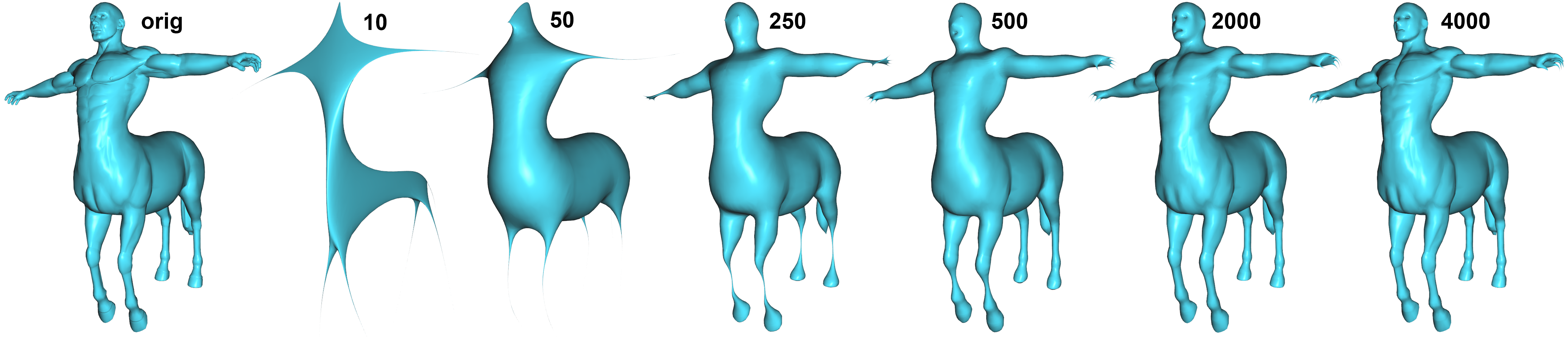}
    \caption{Geometric reconstruction of the Centaur model (left-most) using an increasing number of Laplace--Beltrami eigenfunctions. A sufficient number of eigenfunctions is required to obtain reconstruction that preserves details of the shape.}
    \label{fig:projection} 
\end{figure*}

\subsection{Lanczos and preconditioned eigensolver}
In Section 6 of \cite{Nasikun2021}, the timings of HSIM are compared with the timings of Lanczos solvers and LOBPCG. Table~\ref{tab:comparison} shows additional results that complement Table~2 in \cite{Nasikun2021}.
\subsection{Fast Approximation}
We compare HSIM with the fast approximation method from \cite{Nasikun2018}. The approximation method has the advantage that the computation times are much shorter and storing the approximate eigenfunctions requires less memory. On the other hand, the approximation errors of \cite{Nasikun2018} are much larger than the errors resulting from HSIM. 
The top row of Figure~\ref{fig:vsSGP1} shows plots of residuals of eigenpairs computed with the approximation scheme from \cite{Nasikun2018} and compares them with the residuals from HSIM with tolerances $\epsilon=10^{-2}$ and $\epsilon=10^{-4}$. The residuals obtained for the approximation scheme from \cite{Nasikun2018} are $10^0$. In contrast, HSIM allows for controlling the residuals. 
The bottom row of the figure additionally shows the computed eigenvalues.
While visually there is no difference between the two HSIM results, the eigenvalues computed with \cite{Nasikun2018} differ significantly from the results of HSIM. We would like to note that in \cite{Nasikun2018} it is advised to use only the first half of the computed eigenvalues. However, significant deviations can be observed in the first half as well.


\begin{figure}[b]
    \includegraphics[width=0.89\linewidth]{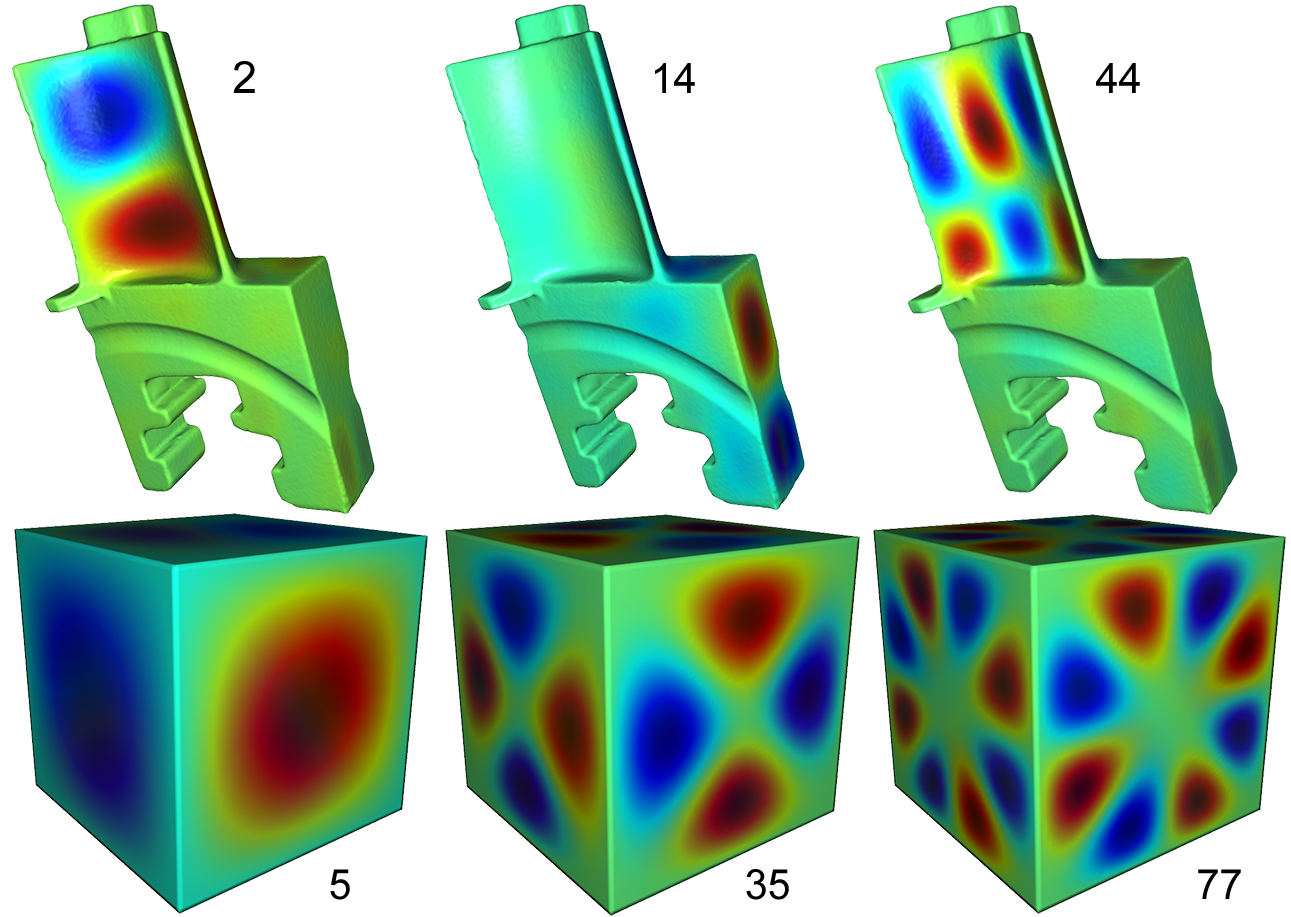}
     \caption{Eigenfunctions of the Hamilton operator are shown.}
     \label{fig:hamiltonian} 
\end{figure}

\begin{figure}[b]
    \subfloat[Gargoyle\label{fig:HKS-Gargoyle-10}]{
        \includegraphics[width=0.47\linewidth]{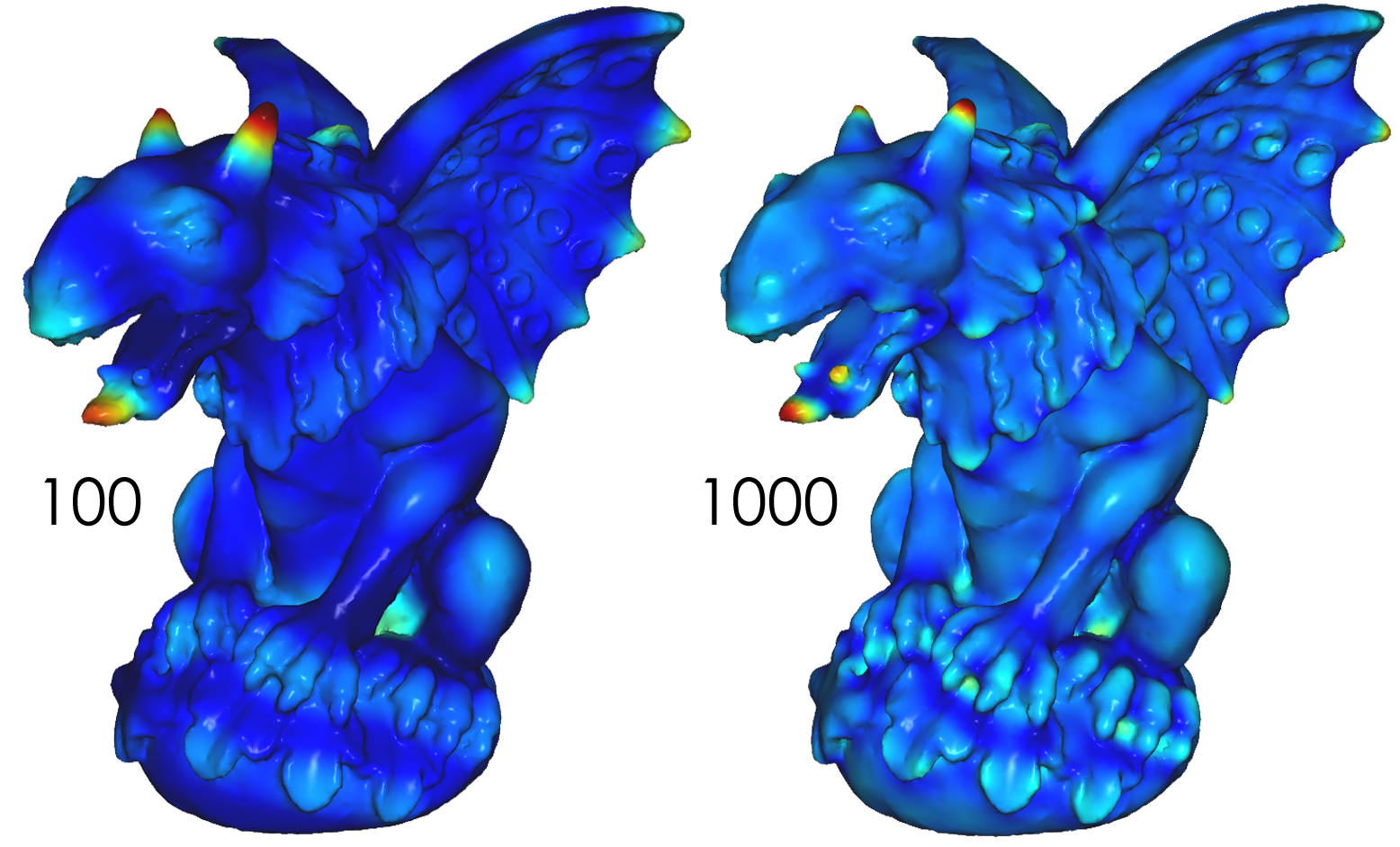}
    }
    \subfloat[Chinese Dragon\label{fig:HKS-Gargoyle-25}]{
        \includegraphics[width=0.49\linewidth]{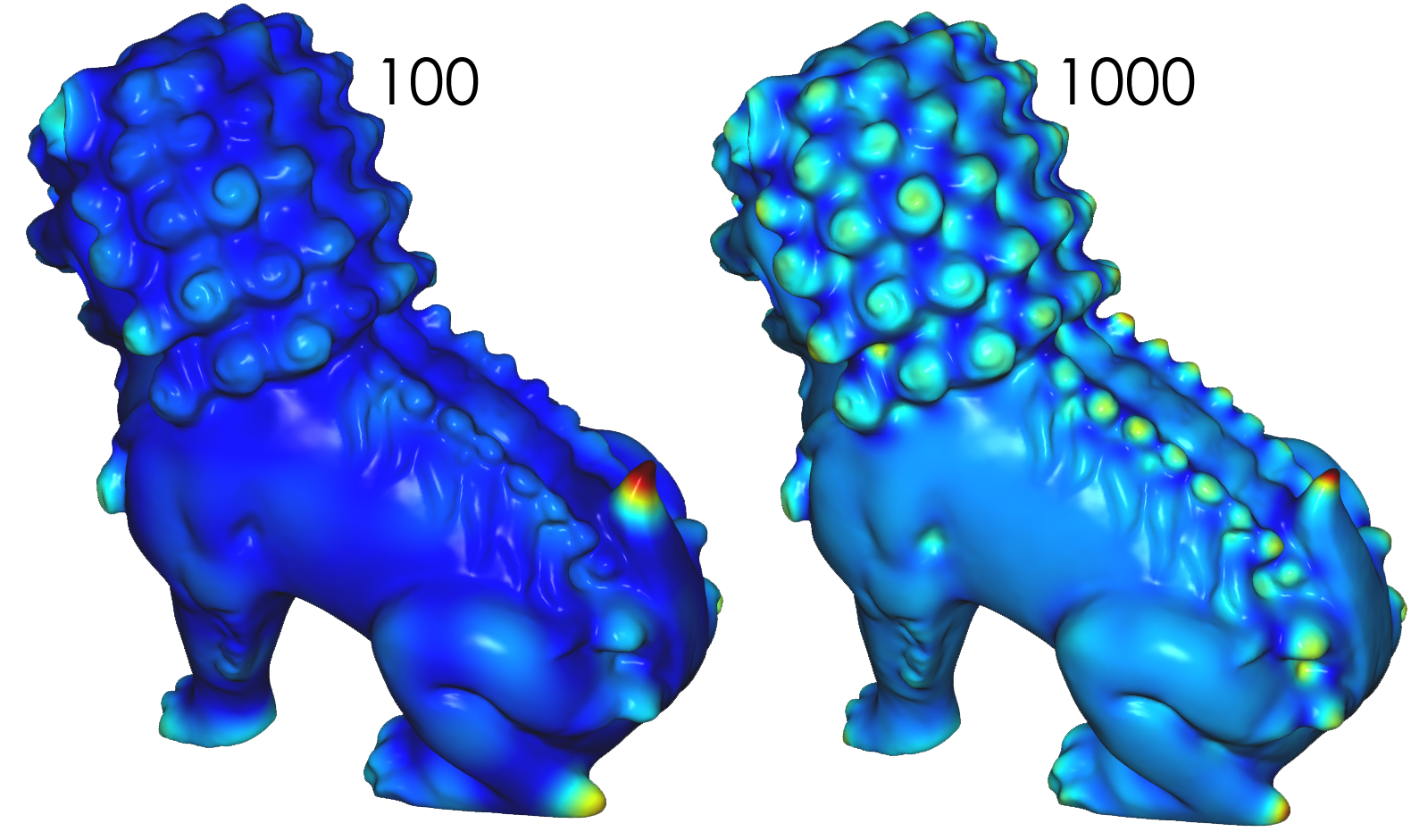}
    }
    \caption{Heat Kernel Signatures (HKS) computed with 100 and 1000 eigenpairs are shown.} 
    \label{fig:HKS-1} 
\end{figure}

\section{Generalization}
\subsection{Hamiltonian operators}
Our evaluation of HSIM is focused on Laplace--Beltrami eigenproblems. In this section, we consider a related operator, the Hamiltonian operator, and present some results for solving Hamiltonian eigenproblems using HSIM. For a background on Hamiltonian operators and their use in spectral analysis, we refer to \cite{Choukroun2020}. 
The Hamilton operators on surfaces we consider are of the form
\[
H:u\rightarrow\Delta u+Vu,
\]
where $\Delta$ is the Laplace--Beltrami operator and $V$ a scalar potential
function. For our experiments, we used the scalar potential
\[
V=t(\kappa_{1}^{2}+\kappa_{2}^{2}),
\]
where $t\in\mathbb{R}^{\geq0}$ and $\kappa_{1}$ and $\kappa_{2}$ are the principal curvatures. The
eigenmodes of this operator have been studied in the context of shape analysis
in \cite{Hildebrandt2010,Hildebrandt2012}. In contrast to the Laplace--Beltrami eigenfunctions,
the eigenfunctions of this operator depend not only on the intrinsic properties of
the surface but also on its extrinsic curvatures. 
Even for $t=0.1$, the eigenfunctions of this operator are fundamentally 
different from those of the Laplace--Beltrami operator as illustrated in 
Figure \ref{fig:hamiltonian}. Table
\ref{tab:hamiltonian} lists iteration counts and timings for solving 
Hamiltonian eigenproblems for $\alpha=0.1$ and $\alpha=1$. As a reference, we also
list the timings for the corresponding Laplace--Beltrami eigenproblem. For $\alpha=0.1$, we
obtain almost the same timings as for the Laplace--Beltrami eigenproblems and for $\alpha=1$, we
noticed in some cases an increase of the required computation time of up to
30\%.

\section{Applications}
In this section, we consider methods that use the Laplace--Beltrami eigenfunctions for shape analysis and processing. We demonstrate that the methods can benefit from using a larger number of eigenfunctions. HSIM facilitates the computation of larger numbers of eigenfunctions.

\subsection{Shape Signatures}
We first consider the Heat Kernel Signature~\cite{SOG09} as an example of a shape signature.
Figure~\ref{fig:HKS-1} shows the Heat Kernel Signature color-coded on two meshes. For both meshes, results using 100 and 1000 eigenfunctions are shown. One can see that the surface details such as the curls of the Chinese lion model are better resolved when 1000 eigenfunctions are used. 
As a consequence, the Heat Kernel Distance delivers better results for finding similar points on a surface when more eigenfunctions are used. Figure~\ref{fig:HKS-region} shows results where similar points to a given point are searched. The results are shown by binary color-coding, where similar points are orange. On the Armadillo mesh, a point at a fingertip is given and on the dinosaur mesh, a point on a toe is given. It can be seen that if 1000 eigenfunctions are used, on both meshes all fingertips and toe tips are found. For 100 eigenfunctions this is not the case. Only about half of the fingertips and toes are found. 
\subsection{Projection}
Methods such as mesh filtering~\cite{Vallet2008} and mesh compression~\cite{Karni2000} need to project the embedding of a surface to the space spanned by the lowest $n$ eigenfunctions. Figure~\ref{fig:projection} shows the results of this projection for the centaur mesh with different values of $n$ ranging from 10 to 4000. One can see that the higher the number of eigenfunctions is, the more surface details are preserved. Even when 2000 eigenfunctions are used, the resulting projection is visually smoother than the original mesh.


\bibliographystyle{ACM-Reference-Format}
\bibliography{Multilevel-Eigensolver}

\end{document}